\documentclass[preprint]{elsarticle}
\usepackage{color} 
\usepackage{hyperref}

\journal{arXiv}

\usepackage{amsmath,amssymb}
\usepackage{graphicx}
 \usepackage{yhmath}

\renewcommand{\theequation}{\arabic{section}.\arabic{equation}}

\input{diagrams.sty}

\renewcommand{\Bar}{\overline}
\newcommand{\bfr}{{\bf R}}
\newcommand{\dotprod}{{\mbox{
\hspace{-.06in}\raisebox{-.015in}{\bf {\Large $\cdot$}}}}}

\newcommand{\id}{{\rm id.}}

\newcommand{\intersect}{\mbox{\small \ $\bigcap$\ }}
\newcommand{\iso}{\cong}
\newcommand{\lb}{\langle}
\newcommand{\rb}{\rangle}

\newcommand{\plus}{\oplus}

\newcommand{\sgn}{{\rm sgn}}
\newcommand{\spn}{{\rm span}}

\renewcommand{\Tilde}{\widetilde}

\newcommand{\tr}{{\rm tr}}
\newcommand{\union}{\mbox{\small $\bigcup$}}
\newcommand{\Union}{\bigcup}
\newcommand{\be}{\begin{equation}}
\newcommand{\ee}{\end{equation}}
\newcommand{\bearray}{\begin{eqnarray}}
\newcommand{\eearray}{\end{eqnarray}}
\newcommand{\bestar}{\begin{eqnarray*}}
\newcommand{\eestar}{\end{eqnarray*}}
\newcommand{\ben}{\begin{displaymath}}
\newcommand{\een}{\end{displaymath}}

\newtheorem{theorem}{Theorem}[section]
\newtheorem{thm}[theorem]{Theorem}

\newtheorem{prop}[theorem]{Proposition}
\newtheorem{lemma}[theorem]{Lemma}

\newtheorem{cor}[theorem]{Corollary}
\newtheorem{defn}[theorem]{Definition}
\newtheorem{conjecture}[theorem]{Conjecture}
\newtheorem{remark}[theorem]{Remark}
\newtheorem{remarks}[theorem]{Remarks}

\newtheorem{notation}[theorem]{Notation}

\newtheorem{Example}[theorem]{Example}

\newtheorem{question}[theorem]{Question}

{\addtocounter{theorem}{1}%
\noindent {\bf
Example \arabic{section}.\arabic{theorem} }}{\\[2 ex]}

\def\qedns{\mbox{\hspace{1ex}}
\mbox{\hfill\vrule height10pt width10pt} \\ \vspace{0.5 ex}
}

\newcommand{\pf}{\noi{\bf Proof}: }

\renewcommand{\b}{\beta}
\renewcommand{\d}{\delta}
\newcommand{\e}{\epsilon}
\newcommand{\g}{\gamma}
\renewcommand{\i}{\iota}
\renewcommand{\l}{\lambda}
\renewcommand{\th}{\theta}

\newcommand{\noi}{\noindent}
\renewcommand{\ss}{\vspace{.1in}}

\newcommand{\ssn}{\vspace{.1in}\noindent}

\newcommand{\bfc}{{\bf C}}
\newcommand{\Sop}{SO(p)}
\newcommand{\So}{SO}
\newcommand{\sop}{{\mathfrak {so}}(p)}
\newcommand{\so}{{\mathfrak {so}}}
\newcommand{\frakg}{{\mathfrak g}}

\newcommand{\ad}{{\rm ad}}
\newcommand{\sympp}{{\rm Sym}^+(p)}

\newcommand{\ts}{{\tilde{S}}}
\newcommand{\tsp}{\ts_p}

\newcommand{\D}{{\rm Diag}}

\newcommand{\Dpp}{{\rm Diag}^+(p)}
\newcommand{\E}{{\cal E}}

\newcommand{\Gr}{{\rm Gr}}
\newcommand{\I}{{\cal I}}
\newcommand{\J}{{\sf J}}
\newcommand{\K}{{\sf K}}
\renewcommand{\L}{{\Lambda}}
\newcommand{\M}{{\cal M}}

\newcommand{\Q}{{\cal Q}}
\newcommand{\sfr}{{\sf R}}
\renewcommand{\S}{{\cal S}}
\newcommand{\SR}{{\cal SR}}

\newcommand{\dsr}{d_{\cal SR}}
\newcommand{\partpset}{{\rm Part}(\{1,\dots,p\})}
\newcommand{\partp}{{\rm Part}(p)}
\newcommand{\Stab}{{\rm Stab}}

\newcommand{\boldsig}{{\mbox{\boldmath $\sigma$ \unboldmath}
           \mbox{\hspace{-.05in}}}}
\newcommand{\miniboldsig}{{\mbox{\scriptsize \boldmath $\sigma$ \unboldmath}
           \mbox{\hspace{-.05in}}}}
\newcommand{\mbs}{\miniboldsig}
\newcommand{\ztwo}{{\bf Z}_2}

\newcommand{\bfe}{{\bf e}}

\newcommand{\bfv}{{\bf v}}
\newcommand{\bfw}{{\bf w}}

\newcommand{\diag}{{\rm diag}}
\newcommand{\bone}{{\bf 1}}
\newcommand{\incl}{{\rm incl}}
\newcommand{\proj}{{\rm proj}}
\newcommand{\sproj}{\mbox{\small ${\rm proj}_2$}}
\newcommand{\bsproj}{\mbox{\small $\Bar{{\rm proj}_2}$}}

\newcommand{\Comp}{{\rm Comp}}
\newcommand{\cald}{{\cal D}}
\newcommand{\gsop}{g_{\So}}
\newcommand{\gdp}{{g_{\cald^+}}}
\newcommand{\dsop}{d_{\So}}
\newcommand{\dgr}{d_{Gr}}
\newcommand{\ddp}{d_{\cald^+}}

\newcommand{\jtop}{{\J_{\rm top}}}

\newcommand{\lbl}{{\rm lbl}}
\newcommand{\slbl}{\mbox{\small ${\rm lbl}$}}
\newcommand{\bslbl}{\mbox{\small $\Bar{{\rm lbl}}$}}
\newcommand{\quo}{{\rm quo}}
\newcommand{\rsym}{R_{\rm sym}}
\newcommand{\jmp}{{\cal J}_{m,p}}
\newcommand{\floor}[1]{\lfloor #1 \rfloor}
\newcommand{\ceil}[1]{\lceil #1 \rceil}
\newcommand{\lev}{{\rm level}}

\begin{document}

\begin{frontmatter}
\title{ Uniqueness questions in
a scaling-rotation geometry  on
 the space of symmetric positive-definite matrices
}
\tnotetext[mytitlenote]{This work was supported by NIH grant
  R21EB012177 and NSF grant DMS-1307178.}

\author[ufl]{David Groisser\corref{cor1}}
\ead{groisser@ufl.edu}
\author[pitt]{Sungkyu Jung}
\ead{sungkyu@pitt.edu}
\author[ucsd]{Armin Schwartzman}
\ead{armins@ucsd.edu}
\cortext[cor1]{Corresponding author}

\address[ufl]{Department of Mathematics, University of Florida, Gainesville, FL 32611, USA}
\address[pitt]{Department of Statistics, University of Pittsburgh, Pittsburgh, PA 15260, USA}
\address[ucsd]{Division of Biotatistics, University of California, San Diego, CA 92903, USA }

\begin{abstract}
 Jung et al. (2015) introduced
a geometric structure on  ${\rm Sym}^+(p)$, the set of $p \times p$ symmetric
positive-definite
matrices, based on eigen-decomposition.  Eigenstructure determines both a stratification of ${\rm Sym}^+(p)$, defined by eigenvalue multiplicities, and fibers
of the ``eigen-composition" map  $F:M(p):=SO(p)\times{\rm Diag}^+(p)\to{\rm Sym}^+(p)$.
When $M(p)$ is equipped with a suitable Riemannian metric,
the fiber structure
 leads to
notions of {\em scaling-rotation distance} between $X,Y\in \sympp$,
the distance in $M(p)$ between fibers $F^{-1}(X)$ and $F^{-1}(Y)$, and
{\em minimal smooth scaling-rotation (MSSR) curves},
images in ${\rm Sym}^+(p)$ of minimal-length geodesics connecting two fibers.
In this paper we study the geometry of the triple $(M(p),F,\sympp)$,
focusing on some basic questions: For which $X,Y$ is there a {\em unique} MSSR  curve from $X$ to $Y$?
More generally, what is the set ${\cal M}(X,Y)$ of
MSSR curves from $X$ to $Y$?
This set is influenced by two
potential types of
non-uniqueness.
We translate the question of whether
the second type
 can occur  into a question about the geometry of Grassmannians $G_m(\bfr^p)$,  with $m$ even,
that we answer for
$p\leq 4$ and $p\geq 11$.
Our method of proof also yields an interesting half-angle formula concerning
principal angles between  subspaces of $\bfr^p$ whose dimensions
may or may not be
equal.
The general-$p$ results concerning MSSR curves and scaling-rotation distance
that we establish here
underpin the explicit
$p=3$ results in
Groisser et al. (2017).
Addressing the uniqueness-related questions requires
a thorough understanding of the fiber structure of $M(p)$, which we also provide.

\end{abstract}


\end{frontmatter}

\section{Introduction}

In this work, we investigate a geometric structure on  $\sympp$, the set of $p \times p$ symmetric
positive-definite (SPD) matrices, $p>1$, and special curves that this structure gives rise to.
Both the geometric structure and these special curves are built from eigen-decomposition of SPD
matrices.

Let $\D^+(p)$ denote the set of $p\times p$ diagonal matrices with
positive diagonal entries. By an {\em (orthonormal) eigen-decomposition} of $X\in\sympp$ we will mean a pair
$(U,D)\in\Sop\times\D^+(p)$ such that $X=UDU^{-1}=UDU^T$.  The space of such
decompositions,

\be
M(p):=\Sop\times \D^+(p),
\ee

\noi  thus comes naturally equipped with a smooth
surjective map $F: M\to\sympp$ defined by

\be
F(U,D) = UDU^T.
\ee

For each $X\in\sympp$ we call
the set $\E_X := F^{-1}(X)$  the {\em fiber over $X$}.  However,  $M(p)$ is not a fiber bundle over $\sympp$ with projection $F$;
the map $F$ is
not even a submersion.  (Rather, the relation of $M(p)$ to $\sympp$ is reminiscent
of the notion of {\em blow-up} in algebraic geometry: $M(p)$ can be viewed as a sort of blow-up
of $\sympp$ along several subvarieties.)
The natural action $\Sop\times \sympp\to\sympp$, $(U,X)\mapsto UXU^T$,
endows $\sympp$ with a stratification by orbit-type, and the derivative of $F$ is nonsingular
only on the pre-image of the``top" stratum. This stratification is identical to the
stratification by ``eigenvalue-multiplicity type", in which the strata are labeled by
partition of the integer $p$. Eigenvalue multiplicities also determine a more refined
stratification of the space $M(p)$, in which the strata are labeled by partitions
of the {\em set} $\{1,\dots,p\}$.  Appendix B
reviews these stratifications.

The fiber structure  of $M(p)$ formalizes the notion of {\em minimal smooth scaling-rotation  curves} \cite{JSG2015scarot}. In 2006, motivated by applications to diffusion-tensor imaging,
Schwartzman \cite{Schwartzman2006} introduced {\em smooth scaling-rotation
curves} as a way of interpolating between SPD matrices in such a way
that eigenvectors and eigenvalues both change at uniform speed.
{\em Minimal} smooth scaling-rotation  curves were defined in
\cite{JSG2015scarot} as  smooth curves of shortest length as determined by an
appropriate Riemannian metric on $M(p)$---curves that
minimize a suitable measure of the
amount of scaling and rotation needed to transform one SPD matrix into another.

More precisely, each factor of $M(p)$ is a Lie group, and for our Riemannian metric $g_M$
on $M(p)$
we take a product metric determined by choosing bi-invariant metrics
$\gsop, \gdp$ on the factors.
We define {\em smooth scaling-rotation (SSR) curves} in $\sympp$
to be the projections to $\sympp$ of geodesics in $(M(p),g_M)$.
In this scaling-rotation framework, the ``distance''
$\dsr(X,Y)$ between any two matrices  $X,Y\in \sympp$  is
defined to be the distance between the fibers $\E_X$ and $\E_Y$ (nonzero if $X\neq Y$ since
each fiber is compact).   We use the term {\em $F$-minimal geodesic} for
a minimal-length geodesic connecting two fibers
$\E_X$ and $\E_Y$,  and {\em minimal pair} for the pair of endpoints of such
a geodesic. A {\em minimal
  smooth scaling-rotation}  (MSSR) curve is the image under $F$ of an $F$-minimal geodesic.

As shown in \cite{JSG2015scarot},
$\dsr$ restricts to a metric on the top stratum of $\sympp$, but is not a metric
on all of $\sympp$.  In \cite{GJS2016metric}, we show that $\dsr$ generates a true metric
$\rho_{\cal SR}$ on $\sympp$ and investigate features of this metric. But
fully understanding the geometry of the
metric $\rho_{\cal SR}$ relies on first understanding MSSR curves,
the function $\dsr$, and related issues we address in the present article.

This paper is devoted primarily to uniqueness-related issues that arise in studying MSSR curves,
and to some unanticipated geometric
results (described in more detail below), potentially of independent
interest, that   were discovered as a result of studying these issues.

 A thorough understanding of the fibers of $F$
is key
to analyzing
several features of
the scaling-rotation framework, including these uniqueness-related issues.
Appendix A provides a thorough picture of the fiber structure of $M(p)$,
including its inextricable tie to the group $\tsp^+$ of ``even
signed-permutations'',  a group not to be confused with a more familiar
group of the same
order and similar-sounding description in terms of signs and permutations, the Weyl group of the simple Lie
algebra $D_p$.
Some results proven
in Appendix A
are applied earlier in the main body of this paper, and some were previously stated
without proof in \cite{GJS2017eigen_EJS} and applied there.

The uniqueness-related results in this paper contribute to
a rigorous and systematic description of the geometry and topology
of the triple \linebreak  $(M(p),F,\sympp)$,
and to a firm foundation for further study of the scaling-rotation
framework, such as in \cite{GJS2017eigen_EJS} and \cite{GJS2016metric}.

 Some of the uniqueness issues we study are related
directly to (non-)uniqueness of MSSR curves themselves,  while others are related
more directly to \linebreak (non-)uniqueness of minimal pairs.
It is easy to see that for all $X,Y\in \sympp$, at least one MSSR curve from $X$ to
$Y$ exists; however,
such curves are not always unique. The dependence on $X$ and $Y$
 of the set $\M(X,Y)$ of such curves is
quite intricate,   and relates strongly to the stratified nature of $\sympp$.
Non-uniqueness issues for minimal pairs are important because
computing $\dsr(X,Y)$ and MSSR curves from $X$ to $Y$
requires finding minimal pairs in $\E_X\times\E_Y$. Even when
the resulting MSSR curve from $X$ to $Y$ is unique, minimal
pairs in $\E_X\times \E_Y$ are never unique, because $\tsp^+$
acts on $M(p)$ in a nontrivial isometric, fiber-preserving fashion.
This action carries minimal pairs to minimal pairs.  For some $(X,Y)$, there
 are also minimal pairs that are not related to each other by this action.

The broad structure of this paper is as follows. Section \ref{sec:notn_prelims}
establishes notation.
Sections \ref{sec:dsrmssr} and \ref{sect:app} contain the statements of most of
our main results, which we will describe below, and those proofs that can be
given quickly.  The proofs of many of our results---especially the ``bonus'' results
that are applicable
outside the scaling-rotation framework entirely---are quite long; these occupy
Sections \ref{sec:app.3}--\ref{sec:app.5}.

 We devote the remainder of this introduction to a more detailed
outline of the paper, and more detailed descriptions of the questions we study and the
results we achieve.

In Section \ref{sect3.1}
we review the basics of SSR curves, before restricting attention to MSSR curves
in Section \ref{sec:dsrmssr2} and beyond. In Section \ref{sec:dsrmssr2} we
discuss the computational-complexity problem arising from
the non-uniqueness of minimal pairs.  Proposition \ref{distprop} takes
advantage of the $\tilde{S}_p^+$-action  by using double-cosets
in $\tsp^+$ to reduce the complexity of computing $\dsr(X,Y)$, of
characterizing all minimal pairs, and of finding all MSSR curves from $X$ to $Y$.
 Proposition \ref{distprop} was applied in \cite{GJS2017eigen_EJS}
to help derive closed-form formulas for $\dsr$
and MSSR curves for $p=3$.  Also discussed and proven in Section
\ref{sec:dsrmssr2} is a general result about scaling-rotation curves: All
such curves are either constant-maps or immersions. This result
is important for an understanding of MSSR curves.

In Section \ref{sect4.4}, we  begin to address uniqueness questions
for MSSR curves, the most basic of which is: under what conditions
on $X,Y\in\sympp$ is there  more than one MSSR curve
from $X$ to $Y$? A more refined version of this question is: for each
pair $(X,Y)$, what is the set $\M(X,Y)$  explicitly?
By characterizing all minimal pairs, Proposition \ref{distprop} provides
a starting point for answering this question.  Among this proposition's
outcomes is also the fact that,  for each pair $(X,Y)$,
every MSSR curve  from $X$ to $Y$ is represented by a minimal pair
whose first point lies in any given connected component of $\E_X$.
But to completely understand $\M(X,Y)$---or even just determine its
cardinality---we still need a way to tell whether MSSR curves corresponding
to two (not necessarily distinct) minimal pairs with first point in
a given connected component of $\E_X$ are the same.
Proposition \ref{prop:MSRequal} gives a necessary and sufficient criterion.
This result  was applied in \cite{GJS2017eigen_EJS}, where it enabled an
explicit computation of the sets $\M(X,Y)$ for $p=3$ when $X$ and $Y$
do not both lie in the top stratum.

In Section \ref{sect4.4} we also define two different ways that non-uniqueness
of MSSR curves can occur.  Given $X,Y\in \sympp$,
for there to be more than one MSSR curve from $X$ to $Y$,
there must exist distinct shortest-length geodesics
$\g_1,\g_2 : [0,1]\to M(p)$ from $\E_X$ to $\E_Y$ such
that $F\circ \g_1\neq F\circ\g_2$. There are essentially two ways,
not mutually exclusive, that this can happen: (i) there can exist
such $\g_i$ ($i=1,2$) whose endpoint-pairs are distinct  minimal pairs,
and (ii) there exist such $\g_i$ whose endpoint-pairs are the same minimal
pair. We call these possibilities ``Type I" and ``Type II" non-uniqueness,
respectively.  Proposition \ref{prop:MSRequal} applies to both.

The study of Type II non-uniqueness, which we begin in Section  \ref{sect:uniqq},
turns out to be especially fruitful.  A minimal pair $((U,D),(V,\L))\in M(p)\times M(p)$
has more than one minimal geodesic connecting its points if and only if the pair $(U,V)\in\Sop\times\Sop$ is {\em geodesically antipodal} (Definition \ref{geodantip}),
which is equivalent to $V^{-1}U$ being an involution.
Our chief tool for determining whether such minimal pairs exist is a property we call
{\em sign-change reducibilty}:  we say that the pair $(U,V)$ is sign-change reducible if
$\dsop(U,V)$ can be reduced by multiplying $U$ or $V$ by a (positive-determinant) ``sign-change matrix'', a diagonal matrix each of whose diagonal entries is $\pm 1$.

We show in Proposition \ref{nscrminimal} that if $(U,V)\in \Sop\times \Sop$ is
{\em not} sign-change reducible, then there exist $D,\L$ in the top stratum of
$\D^+(p)$ such that $((U,D), (V,\L))$ is a minimal pair.  We show in
Proposition \ref{pleq4prop} that for $p\leq 4$, every geodesically antipodal pair
$(U,V)$ is sign-change reducible, and that for $p\geq 11$, there exist geodesically
antipodal pairs that are not sign-change reducible. From these propositions we
deduce that Type II non-uniqueness never occurs for $p\leq 4$
(Corollary \ref{pleq4cor}), and that it always occurs for {\em some}
$(X,Y)\in\sympp\times \sympp$ if $p\geq 11$ (Corollary \ref{gamp}). We do
not believe that either of the numbers 4 and 11 above is sharp; our methods
are simply not conclusive when $5\leq p\leq 10$.

Together, Proposition \ref{nscrminimal} and Corollary \ref{gamp}
show that sign-change reducibility is the only obstruction to having points
$X,Y$ in the top stratum of $\sympp$ for which the set $\M(X,Y)$
exhibits Type II non-uniqueness.

Even without Proposition \ref{pleq4prop}, for $p\leq 3$  it is rather trivial
that all geodesically antipodal pairs are sign-change reducible, and for $p=4$
an independent proof relying on quaternions is also possible.  However, our proof
of the $p\leq 4$ part of Proposition \ref{pleq4prop} makes no use of quaternions,
and unifies these low-$p$ results.

Our proof of Proposition \ref{pleq4prop}, completed in Section \ref{sec:app.5}
after laying groundwork in Sections \ref{sect:app}--\ref{sec:app.4}, takes us in
unexpected directions, with unanticipated consequences. We initially introduced
the notion of sign-change reducibility into our scaling-rotation-curve study as an
{\em ad hoc} tool to help us determine whether Type II non-uniqueness
of MSSR curves,  impossible for $p\leq 4$, is {\em ever} possible. This is
equivalent to answering the question ``Are all geodesically antipodal pairs
in $\Sop\times\Sop$ sign-change reducible?'' But as we show in Proposition
\ref{equivstat}, a refined version of the  latter question
is equivalent to a question purely about the geometry of Grassmannians
equipped with a standard Riemannian metric: for $m$ even and positive,
is every $m$-dimensional subspace of $\bfr^p$ within a certain
distance $c(m)$ of a coordinate $m$-plane?  (This question can, of course, be
asked without restricting the parity of $m$, but the above equivalence leads us
to consider only even $m$ in this paper.)
By constructing examples, we show that for $m=2$,  the answer to the
Grassmannian question is no for $p\geq 11$. This, combined with the equivalence result
in Proposition \ref{equivstat}, yields the ``$p\geq 11$'' part of Proposition
\ref{pleq4prop} mentioned above.  The ``$p\leq 4$'' part of Proposition
\ref{pleq4prop} is proven by other means
(via the more technical Proposition \ref{involprop-weak}).

  While the possibility of Type-II non-uniqueness is what led us to
the question above  about Grassmannians,  this question and our study of
it may be of independent interest. Our study led us to investigate several related
questions concerning distances between (even-dimensional) subspaces of
$\bfr^p$ and (even-dimensional) coordinate planes not necessarily of the
same dimension.
Perhaps the most unexpected of these  is a half-angle relation stated in Proposition \ref{cor:halfangle} and
proven in Section \ref{sec:app.3}:  for any two involutions $R_1, R_2\in \Sop$, each of the principal angles between the
$(-1)$-eigenspaces of $R_1$ and $R_2$ is exactly half a correspondingly indexed normal-form angle of $R_1R_2$.  This relationship holds whether or not the dimensions of the $(-1)$-eigenspaces are equal.
When the dimensions {\em are} equal, we use this relationship to show that
a natural correspondence between $\Gr_m(\bfr^p)$ and a connected component
of the set of involutions in $\So(p)$
is a metric-space isometry
up to a constant factor of 2 (Proposition \ref{isometry}). This isometric relation is also deducible (and may already be known) from a purely Riemannian approach, but our
proof uses essentially no Riemannian geometry (see Remark \ref{rem:isom} for a more precise
statement, and an additional interpretation of what our proof of Proposition \ref{isometry} shows).

The most important results coming from our study of sign-change reducibility are stated in Section \ref{sect:app}, with the proofs deferred to Sections \ref{sec:app.3}, \ref{sec:app.4}, and \ref{sec:app.5}.  These results include those mentioned above, and one more
whose statement involves terminology not included in this Introduction:
Proposition \ref{conjtruem2}, a special case of a more general conjecture we make about
sign-change reducibility (Conjecture \ref{onlymcanwork}).
Key to almost all of these results
is the technical Lemma \ref{R-eigenlemma},  which establishes several
facts concerning the product of
a general involution in $\Sop$ and a positive-determinant sign-change matrix.

 We mention in passing that there is a vast body of literature
devoted to defining and studying ``distance-functions'' (not necessarily true metrics) on $\sympp$
different from the scaling-rotation distance $\dsr$ and metric $\rho_\SR$; for a discussion and comparison see \cite{GJS2017eigen_EJS}
and the references therein.

\setcounter{equation}{0}
\section{Notational preliminaries}\label{sec:notn_prelims}

In this paper,  when a group $G$ acts from the left on a
set $X$ in a previously specified way, we generally denote the action simply by
$(g,x)\mapsto g\,\dotprod x$.

Let $\partpset$ denote the set of partitions of  $\{1,2,\dots,p\}$, and $\partp$
the set of partitions of the integer $p$. Let $\D(p)$ denote the set of $p\times p$
diagonal matrices. Each $D\in\D(p)$ naturally determines an element
$\J_D\in\partpset$ according to ``which eigenvalues are equal'' (see Notation
\ref{defjd}). The group $\Sop$ acts on $\sympp$ on the left via
$(U,X)\mapsto UXU^T$. The stabilizer $G_D$ of $D\in \D^+(p)$
under this action depends only on $\J_D$, and $G_{D_1}=G_{D_2}$ if
$\J_{D_1}=\J_{D_2}$.  For each $\J\in\partpset$ we may define a subgroup
$G_\J\subset \Sop$ by declaring  to be $G_\J=G_D$ for any $D$ for which
$\J_D=\J$.  We write $G_D^0,G_\J^0$ for the identity component of $G_D,G_\J$ respectively. See Appendix A for an alternative definition of $G_\J$ and additional facts
concerning these groups.

We call an element of $U\in O(p)$ a {\em signed-permutation matrix} if every entry
of $U$ is either $0$ or $\pm 1$, and call a signed-permutation matrix {\em even}
if it lies in $\Sop$.  The set of even signed-permutation matrices
forms a subgroup $\tsp^+\subset \Sop$ of order $2^{p-1}p!$.  As discussed in
Appendix A (Section \ref{gpthy}), we view this subgroup as a canonical copy of
an ``abstract'' group $\tsp^+$ of {\em even signed-permutations}, an extension
of the symmetric group $S_p$.   We will typically denote an even signed-permutation
by the letter $g$, and the corresponding matrix by $P_g$.  We denote the
natural epimomorphism $\tsp^+\to S_p$ by $g\mapsto \pi_g$.  The group
$\tsp^+$ plays a critical role in understanding the fibers of $F$ (starting with Corollary \ref{prop:tsp+_action} in the next section) and in simplifying computations of $\dsr$.
This group, which is not encountered in geometry as often as another
group of the same order, is discussed in greater detail in Appendix A.

We define $\I_p^+=\tsp^+\intersect \D^+(p)$, and call elements of $\I_p^+$
 {\em (even) sign-change matrices}.   We view $\I_p^+$ as
a copy of a (certain) index-two subgroup of $({\bf Z}_2)^p$, as discussed in
Appendix B.  We will typically denote an element of the abstract group $\I_p^+$ by the letter $\boldsig$, and the corresponding matrix by $I_\mbs$.

\begin{notation}\label{defkj}
\rm (a)
For $\J=\{J_1,\dots, J_r\}\in\partpset$, define  (i) $\Gamma_\J = \tsp^+\intersect G_\J$,
(ii) $\Gamma_\J^0 = \Gamma_\J\intersect G_\J^0 = \tsp^+\intersect
G_\J^0$, (iii) $K_\J = \{\pi\in S_p :\pi(J_i)=J_i, \
\linebreak 1\leq i\leq r\}\subset S_p$,
and (iv) $\tilde{K}_\J=\{g\in \tsp^+: \pi_g\in
K_\J\}$. Observe that $\I_p^+\subset\tilde{K}_\J$, and that $K_\J
=\{\pi\in S_p\mid \pi\dotprod D =D \ \mbox{for {\em some} $D$ with $\J_D=\J$}\}=\{\pi\in S_p\mid \pi\dotprod D =D \ \mbox{for {\em all} $D$ with $\J_D=\J$}\}.$

(b)  For any $X\in\sympp$ and $(U,D)\in \E_X$, define
\be\label{defcompud}
[(U,D)]=\{(UR,D) : R\in G_D^0\},
\ee
\noi the connected component of $\E_X$ containing $(U,D)$. We write
$\Comp(\E_X)$ for the set of connected components of $\E_X$.

(c) For any Lie group $G$ and closed subgroup $K$, we write $G/K$ and
$K\backslash G$ for the spaces of left- and right-cosets,
respectively, of $K$ in $G$.
\end{notation}

\setcounter{equation}{0}
\section{The scaling-rotation framework and some results
for scaling-rotation curves}
\label{sec:dsrmssr}

The Lie groups $\So(p)$ and $\D^+(p)$ carry natural bi-invariant Riemannian
metrics.  If we endow  $M(p)=\So(p)\times\D^+(p)$ with a product Riemannian
metric $g_M$  constructed from these, the geodesics $\g$ in $(M(p),g_M)$ are
easily computed. We define {\em smooth scaling-rotation (SSR) curves} in $\sympp$
to be the projections to $\sympp$ of the geodesics in $(M(p),g_M)$, i.e. curves of
the form $F\circ \g$. (In \cite{Schwartzman2006} and \cite{JSG2015scarot} these
were called simply ``scaling-rotation curves".
In Section \ref{sec:dsrmssr2} we explain why we have
added ``smooth" to this name.)

\subsection{Smooth scaling-rotation curves}\label{sect3.1}

The Lie algebra $\sop=T_I(\Sop)$ is the space of $p\times p$
antisymmetric matrices.
The bi-invariant Riemannian metric
$\gsop$ on $\Sop$ we will use is defined at the identity $I\in\Sop$ by
\be\label{gsop}
\left.\gsop\right|_I(A_1,A_2)=-\frac{1}{2}\tr(A_1 A_2),
\ee
(The requirement of bi-invariance determines a Riemannian metric on $\Sop$
up to a constant factor unless $p=4$,  of course, but for all $p\geq 3$ the inner
product \eqref{gsop} is a multiple of the Killing form.)

Since the abelian Lie group $\D^+(p)$ is an open subset of the vector space
$\D(p)$, for each $D\in\D^+(p)$ we will identify $T_D(\D^+(p))$ canonically
with $\D(p)$ . With this  identification understood,  the invariant Riemannian
metric $\gdp$ we use is defined by
\be\label{gdp}
\left.\gdp\right|_D(L_1,L_2)=\tr(D^{-1}L_1 D^{-1}L_2)
\ee
\noi where $D\in \D^+(p)$ and $L_1,L_2\in T_D(\D^+(p)).$  Up to a constant factor, $\gdp$
is the unique (bi-)invariant metric on $\D^+(p)$ that is also invariant
under the natural action of the symmetric group $S_p$.

Naturally identifying of $T_{(U,D)}M(p)$ with $T_U(\Sop)\oplus T_D(\D^+(p))$,
the Riemannian metric $M(p)$  we will use is
\be\label{gm}
g_M:=k\gsop\oplus\gdp\ ,
\ee
where
$k>0$ is an arbitrary parameter that can be chosen as desired for applications.

\begin{defn}\label{defssr}
 {\rm A {\em smooth scaling-rotation (SSR) curve} is a
    curve $\chi$ in $\sympp$ of the form $F\circ\g$, where $\g:I\to M(p)$
    is a geodesic defined on some interval $I$. }
\end{defn}

In this paper, we use {\em curve} sometimes to mean a {\em parametrized curve}
(a map with domain some interval), and sometimes to mean an equivalence class of
such maps, where two maps are regarded as equivalent if one is a monotone
reparametrization of the other.  Also, we use the noun {\em geodesic}
sometimes to mean a complete geodesic and sometimes to mean a geodesic segment.
Our intended meanings should always be clear from context.

The geodesics $\g$  in $M(p)$ are exactly the curves of the form
$t\mapsto (\g_1(t),\g_2(t))$, where $\g_1$ is a geodesic in $(\Sop,\gsop)$
and $\g_2$ is a geodesic in $(\D^+(p),\gdp)$. Since the metrics $\gsop$
and $\gdp$ are bi-invariant, the geodesics in $(\Sop,\gsop)$ and
$(\D^+(p),\gdp)$ can be obtained as either left-translates or right-translates
of geodesics through the identity. For agreement with \cite{JSG2015scarot} and \cite{GJS2017eigen_EJS}, in this paper we use right-translates.

It well known that in the Riemannian manifold $(\Sop, \gsop)$, the
cut-locus of the identity is the set of all involutions, $\{R\in\Sop\mid R^2=I\neq R\}$.
For every non-involution $R\in\Sop$, there is a unique $A\in\sop$ of
smallest norm such that $\exp(A)=R$ (see Section \ref{sec:app.1});
we define $\log(R)=A$. If $R$ is an involution, there  is not a unique such
$A$, but all minimal-norm $A$'s with $\exp(A)=R$ have the same
norm, which we denote $\|\log(R)\|$.  (Thus $\|\log(R)\|$ is a
well-defined real number for all $R\in\Sop$, even
when there is no uniquely defined element``$\log R$'' in $\sop$.)
With this understood, the geodesic-distance function $d_M$ on $M(p)$ is given by
\bearray
\label{defdM}
d_M^2\left((U,D),(V,\Lambda)  \right)
& = & k\, \dsop(U,V)^2 + \ddp(D, \Lambda)^2
\\
&=& \frac{k}{2}\left\|\log (U^{-1}V)\right\|^2 +
\left\|\log( D^{-1}\L)\right\|^2,
\label{dM2}
\eearray
\noi where in \eqref{dM2} and for the rest of this paper, $\| \ \|$
denotes the Frobenius norm on matrices: $\|A\|^2=\|A\|_F^2=\tr(A^TA)$
for any matrix $A$.

The invariances of the metrics $\dsop$ and $\ddp$ lead to
the following proposition, key to many of our results
(e.g. Proposition, \ref{distprop}, Proposition \ref{prop:biject}, and
Corollary \ref{fibdescrip}).

\begin{prop}\label{prop:tsp+_action}
The map $\tsp^+\times M(p)\to M(p)$ defined by
\be\label{act1}
(g,(U,D)) \mapsto g\dotprod(U,D) :=(UP_g^{-1}, \pi_g\dotprod D)
\ee
is a free, isometric, left-action of $\tsp^+$ on $M(p)$ that preserves every fiber of $F$. \qedns
\end{prop}

\subsection{
Scaling-rotation distance and MSSR curves}\label{sec:dsrmssr2}

\begin{defn}[{\cite[Definition 3.10]{JSG2015scarot}}]{\rm For $X,Y\in \sympp$,
        the {\em scaling-rotation distance} $\dsr(X,Y)$ between $X$
        and $Y$ is defined by

\be\label{defdsr} \dsr(X,Y) := \inf_{ \substack{(U,D) \in \E_X,
    \\ (V,\Lambda) \in \E_Y }} d_M( (U,D), (V,\L) ).
\ee
}
\end{defn}

\begin{defn} \label{defmssr}{\rm Let $\g$ be a piecewise-smooth curve in $M(p)$ and let
    $\ell(\g)$ denote the length of $\g$.  For $X,Y\in\sympp$, we call
    $\g:[0,1]\to M(p)$ an {\em $F$-minimal geodesic (from $\E_X$ to $\E_Y$)} if
    $\g(0)\in \E_X, \g(1)\in \E_Y$, and $\ell(\g)=\dsr(X,Y)$.  We call
    a pair of points $((U,D),(V,\L))\in \E_X\times \E_Y$ a {\em
      minimal pair} if $(U,D)=\g(0)$ and $(V,\L)=\g(1)$ for some
    $F$-minimal geodesic $\g$. A {\em minimal smooth
      scaling-rotation} (MSSR) curve from $X$ to $Y$ is a curve $\chi$
    in $\sympp$ of the form $F\circ\g$ where $\g$ is an
    $F$-minimal geodesic. We say that the MSSR curve $\chi=F\circ\g$
{\em corresponds to} the minimal pair formed by the endpoints of $\g$.
We let $\M(X,Y)$ denote the set of MSSR curves from $X$ to $Y.$}
\end{defn}

Obviously an $F$-minimal geodesic is a minimal geodesic in
the usual sense of Riemannian geometry: it is a curve of shortest
length among all piecewise-smooth curves with the same endpoints. (From
the general theory of geodesics, the image of any such curve $\g$ is actually
{\em smooth}.)
Thus a definition equivalent
to \eqref{defdsr} is

\be\label{defdsr2}
\dsr(X,Y) = \inf\left\{\ell(\g) \mid \g:[0,1]\to
M(p)\ \mbox{is a geodesic with}\  \g(0)\in\E_X\ ,\ \g(1)\in
\E_Y\right\}.
\ee

\noi Thus an $F$-minimal geodesic can alternatively be
defined as a geodesic of minimal length among all geodesics starting
in
one given fiber and ending in another.

Every fiber of $F$ is compact (an explicit description
is given in Corollary \ref{fibdescrip}), so
the infimum in \eqref{defdsr} is always achieved. Hence for all
$X,Y\in \sympp$, there always exists an
$F$-minimal geodesic, a minimal
pair in $\E_X\times \E_Y$, and an MSSR curve from $X$ to $Y$.

\begin{remark}\label{lengthchi}
{\rm Observe that we have not defined a Riemannian metric on \linebreak $\sympp$, so there is no ``automatic" meaning attached to the phrase {\em length of a smooth
curve} in $\sympp$.   However, for an SSR  curve $\chi$ in $\sympp$ we define the {\em length of $\chi$} to be $\ell(\chi):=\inf\{\ell(\g): \g \ \mbox{is a geodesic in $M(p)$ and}\ F\circ\g =\chi\}$.
With this definition, \eqref{defdsr2} becomes
\bearray\nonumber
\dsr(X,Y) &=& \inf\left\{\ell(\chi) \mid \chi:[0,1]\to
\sympp\ \mbox{is an SSR curve with}\right.\\
&&  \left.\phantom{\inf\left\{\ell(\chi) \mid \right.}
\  \chi(0)=X,\,\chi(1)=Y\right\}.
\label{defdsr3}
\eearray
}
\end{remark}

{\em A priori}, given $X,Y\in \sympp$, a concrete
computation of $\dsr(X,Y)$ involves computing the distance in $M(p)$ between each connected
component of $\E_X$ and each connected component of $\E_Y$, then taking the
minimum over all component-pairs.  For $X=F(U,D)$, the number of connected components of $\E_X$ is
$|\tsp^+|/|\Gamma^0_{\J_D}|$ (see Proposition \ref{prop:biject} in Appendix A), which tends
to be a rather large number (see Corollary \ref{fibdescrip}).
It is obvious from Propositions \ref{prop:tsp+_action} and  \ref{prop:biject}) that computing {\em all} the distances
between fiber-components is redundant. It is not so obvious exactly how much redundancy
there is (more than one might guess just from looking at these two propositions).
   As a practical matter, it is desirable
to reduce the number of component-pair computations as much as possible,
taking advantage of less-obvious redundancy.
We will do this in Proposition \ref{distprop} below.
This proposition plays a crucial role in \cite{GJS2017eigen_EJS}, where for $p=3$ we apply it
to compute all scaling-rotation distances, and to
help compute and classify all MSSR
curves.
The proof of Proposition \ref{distprop} (which is given only
in the present paper, not in \cite{GJS2017eigen_EJS})
relies on the characterization
of fibers given in Appendix A as Corollary \ref{cor:fibform4}.

\begin{defn}\label{def:dlbcst}\rm
Recall that given any group $G$ and subgroups $H_1,H_2$, an
$(H_1,H_2)$ {\em double-coset} is an equivalence class under the
equivalence relation $\sim$ on $G$ defined by declaring $g_1\sim g_2$
if there exist $h_1\in H_1, h_2\in H_2$ such that $g_2=h_1g_1h_2$. The
set of equivalence classes under this relation is denoted
$H_1\backslash G/H_2$.  By a {\em set of representatives} of
$H_1\backslash G/H_2$ we mean a subset of $G$ consisting of exactly
one element from each $(H_1,H_2)$ double-coset. Since every left or right
coset is also a double-coset, this defines ``set of representatives''
for ordinary cosets as well.
\end{defn}

\begin{prop}\label{distprop}
Let $X,Y\in\sympp$ and let $(U,D)\in \E_X, (V,\L)\in \E_Y$.  Let $Z$
be any set of representatives of $\Gamma_{\J_D}^0\backslash \tsp^+/
\Gamma_{\J_\L}^0$.  Then the scaling-rotation distance from $X$ to $Y$ is
given by
\begin{align}
\dsr(X,Y)^2  = \min_{g\in Z}
  \left
    \{k\left(
          \wideparen{d}(g ; (U,D) , (V,\Lambda) )
        \right)^2
       +\|\log\left(D^{-1}(\pi_g\dotprod \L)\right)\|^2
      \right\}, \label{fibdist8a}
\end{align}
where
\begin{equation}
\wideparen{d}(g ; (U,D) , (V,\Lambda) )= \min_{R_U\in
                    G_D^0, R_V\in G_\L^0} \left\{ d_{\So}(UR_U, VR_V
P_g^{-1}
) \right\}. \label{fibdist8a-cts-opt}
\end{equation}

Every MSSR curve from $X$ to $Y$ corresponds to some minimal pair
whose first element lies in the connected component $[(U,D)]$ of $\E_X$.
\end{prop}

\pf From Corollary \ref{cor:fibform4}
we have
\be\label{fud}
\E_X\times \E_Y =
\left\{
(g_1\dotprod(UR_U,D), g_2\dotprod(VR_V,\L))
: R_U\in G_D^0, R_V\in G_\L^0; g_1, g_2\in
\tsp^+
\right\}.
\ee
\noi
By Proposition \ref{prop:tsp+_action},
for all
$g_1,g_2\in\tsp^+$ we have
\be
\label{repdist1}
d_M(g_1\dotprod(UR_U,D), g_2\dotprod(VR_V,\L))
=d_M((UR_U,D), (g_1^{-1}g_2)\dotprod(VR_V,\L)).
\ee

Proposition \ref{prop:tsp+_action} implies that the action of $g_1^{-1}$ on $M(p)$
carries a geodesic $\g_1$ with endpoints $g_1\dotprod(UR_U,D), g_2\dotprod(VR_V,\L)$
into a geodesic $\g_2$ with endpoints $(UR_U,D), (g_1^{-1}g_2)\dotprod(VR_V,\L)$
and that satisfies $F\circ \g_1=F\circ\g_2$. Hence, every smooth scaling-rotation
(SSR) curve from $X$ to $Y$ is of the form $F\circ\g$ where $\g:[0,1]\to
M(p)$ is a geodesic with $\g(0)=(UR_U,D)\in [(U,D)]$ and $\g(1)\in \E_Y$.

Suppose $\g_1, \g_2$ are two such geodesics, with $\g_i(1)=
(VR_V P_{g_i}^{-1},\pi_{g_i} \dotprod \L)$, $i=1,2$. If $g_2=h_D g_1
h_\L$, with $h_D\in \Gamma_{\J_D}^0$ and $h_\L\in \Gamma_{\J_\L}^0$,
then
\bestar
\lefteqn{
d_M((UR_U, D), g_2\dotprod (VR_V   , \L))
 \ =\ d_M((UR_U,D), h_D\dotprod g_1\dotprod h_\L\dotprod(VR_V,\L))
}
\\
&=&d_M((UR_Uh_D,(\pi_{h_D})^{-1}\dotprod D),
(VR_Vh_\L^{-1} P_{g_1}^{-1},\pi_{g_1}\dotprod \pi_{h_\L}\dotprod\L))
\\
&=&d_M((UR_{U,1},D), (VR_{V,1}\, P_{g_1}^{-1},\pi_{g_1}\dotprod\L))
\eestar
\noi where $R_{U,1}=R_U h_D\in G_D^0$ and $R_{V,1}=R_V h_\L^{-1}\in
G_\L^0$. The same argument as in the preceding paragraph shows that
the SSR curve determined by the pair $((UR_U, D), (VR_V P_{g_2}^{-1},
\pi_{g_2}\dotprod \L))$ is the same as the SSR curve determined
by the pair $((UR_{U,1}, D), (VR_{V,1}\, P_{g_1}^{-1},
\pi_{g_1}\dotprod\L))$. Hence any representative $g\in \tsp^+$ of a given
$(\Gamma_{\J_D}^0,\Gamma_{\J_\L}^0)$ double-coset determines the same set
of SSR curves as does any other representative of that
double-coset. The Proposition now follows. \qedns

We end this subsection with a discussion and results that motivate
our inclusion of the word {\em smooth} in``smooth scaling-rotation curve''.
By its definition, every SSR curve $\chi:I\to \sympp$ is a smooth {\em map},
but it is not clear whether the image of $\chi$ is ``geometrically smooth'', i.e.
locally (in $I$) a smooth submanifold or submanifold-with-boundary of $\sympp$.
For the image of $\chi$ to be geometrically smooth in this sense, $\chi$ must
admit a {\em regular} parametrization, one that is an immersion.  It turns out
that all SSR curves do, except for those whose images are single points:

\begin{prop}\label{prop:immersion}
If $\g$ is a non-constant geodesic, then  $F\circ\g$ is  either an
immersion or a constant map.
\end{prop}

\pf  Let $\g:[0,1]\to M(p)$ be a non-constant $F$-minimal geodesic and let $\chi=F\circ\g$.

Let $(U,D)=\g(0)$ and let $X=\chi(0)=F(U,D)$.  Since $\g$ is a geodesic there exist unique $A\in\sop, L\in \D(p)$ such that $\g(t)=(e^{tA}U, e^{tL}D)$. Non-constancy implies $(A,L)\neq (0,0)$.
Direct computation yields
\ben
\chi'(t)= e^{tA}\left\{ [A, U\L(t)U^T] +UL \L(t)U^T\right\}e^{-tA},
\een
\noi where $\L(t)=e^{tL}D $ and $[\ , \ ]$ denotes matrix commutator.

Suppose that $t_0\in [0,1]$ is such that $\chi'(t_0)=0$. Then
\be\label{mssr-eq2}
[A, U\L(t_0)U^T] +UL \L(t_0)U^T=0.
\ee
Multiplying on left by $U^T$ and on the right by $U$ yields
$
[\tilde{A},\L(t_0)] +L\L(t_0)=0,
$
where $\tilde{A}=U^TAU$.  But because $\L(t_0)$ is diagonal, the diagonal entries of any
commutator $[B,\L(t_0)]$ are zero.  Since $L\L(t_0)$ is a diagonal matrix, this implies that
$[\tilde{A},\L(t_0)] =0 = L\L(t_0)$.   But $\L(t_0)$ is invertible, so the second equality implies $L=0$.
Thus $\L(t)=D$ for all $t$, and plugging this into \eqref{mssr-eq2} with $t=t_0$ we find $[A,X]=0$.
It follows that $X$ commutes with $e^{tA}$ for every $t$. Hence $\chi(t)=e^{tA}UDU^T e^{-tA}=e^{tA}Xe^{-tA}=X$ for all $t$.

Thus either $\chi'(t)$ is nonzero for every $t\in [0,1]$ or $\chi$ is constant.
\qedns

As noted in \cite{JSG2015scarot}, the ``scaling-rotation distance" $\dsr$ is
{\em not} a metric on $\sympp$; it does not satisfy the triangle inequality.
In \cite{GJS2016metric}, we show that the pseudometric $\rho_\SR$
generated by the semimetric $\dsr$ {\em is} a true metric on $\sympp$.
(It is not trivial to show that $\rho_\SR(X,Y)\neq 0$ for $X\neq Y$.)
Effectively, the construction enlarges the class of scaling-rotation (SR)
curves $\chi$ considered in \eqref{defdsr3} from smooth maps to
piecewise-smooth maps (with $\ell(\chi)$ redefined correspondingly).
This definition of the scaling-rotation metric $\rho_\SR$ is analogous to the
definition of  ``distance between two points in a Riemannian manifold": the
infimum of the lengths of {\em piecewise}-smooth curves joining the points.
But some minimal-length SR curves are {\em geometrically} non-smooth
(having corners); an MSSR curve from $X$to $Y$ has minimal length only
among {\em smooth}  scaling-rotation curves from $X$ to $Y$.  (This
phenomenon does not occur in Riemannian geometry; in a Riemannian manifold,
minimal piecewise-smooth curves between two points are  always geometrically
{\em smooth}.) It is for this reason we have made ``smooth'' part of the
terminology used in Definition \ref{defssr}.

\begin{remark}\rm
It seems likely that a non-constant MSSR curve $\chi$ is actually an embedding
(for this, it suffices that $\chi$ be injective, since $[0,1]$ is compact), but we
have not proven this.  There do exist {\em non-minimal} non-constant SSR curves
that are not one-to-one. One example is any nonconstant periodic SSR curve:
$t\mapsto F(\exp(tA)U,D)$ where $(U,D)\in\sympp$ and $A\in\sop$ is
any nonzero element for which there exists $t_1\neq  0$ such that $\exp(t_1A)=I$.
(For $p\leq 3$, the latter condition is redundant.) The restriction of this curve to
$[0,|t_1|\,]$ is an SSR curve of positive length from $(U,D)$ to $(U,D)$.  A
nonperiodic example with $p=2$ is the following. Let $J=\left(
\begin{array}{rr} 0 & -1 \\ 1 & 0 \end{array}\right)$, $U(t)=\exp(t\frac{\pi}{2} J)$,
$D(t)=\left(
\begin{array}{rr} e^{1-t} & 0 \\ 0& e^t \end{array}\right).$ Then the curve
$t\mapsto \g(t):=(U(t),D(t))$ is a geodesic in $M(2)$. Let $\chi$ be the SSR
curve $F\circ\g$.  Then, as the reader may check, if $t_1<t_2$ we have
$\chi(t_1)=\chi(t_2)$ if (and only if) for some integer $n\geq 0$ we have
$t_1=-n$ and $t_2=n+1$.  Now let $n_1,n_2$ be non-negative integers, let
$t_1\in (-n_1-1,-n_1), t_2\in (n_2+1, n_2+2),$ and  let $n=\min\{n_1,n_2\},
X=\chi(t_1)$, and $Y=\chi(t_2)$.  Then $\chi_{[t_1,t_2]}$ is an SSR curve
from $X$ to $Y$ with $n+1$ self-crossings.   Note that the presence of
self-crossings does not directly imply that $\chi_{[t_1,t_2]}$ is not an MSSR
curve:  if we remove the  closed curve $\chi_{[-n,n+1]}$ from $\chi_{[t_1,t_2]}$,
the piecewise-smooth curve $\chi_1$ from $X$ to $Y$ that remains is not an
SSR curve.  (As the reader may check, the set $\{\chi'(-n),
\linebreak\chi'(n+1)\}$ is linearly independent, so $\chi_1$ cannot be
reparametrized as an immersion.  Hence, by Proposition \ref{prop:immersion},
there is no geodesic $\g_1$ in $M(2)$ such that $\chi_1$ can be reparametrized
as $F\circ\g_1$.) Hence $\chi_1$ is not a candidate for an SSR curve from $X$
to $Y$ that is shorter than $\chi$.   However, with a little effort one can check
by direct computation that there is an $F$-minimal geodesic from $X$ to $Y$
that is shorter than $\g|_{[t_1,t_2]}$.  (One can compute the length of the
minimal geodesic from any of the four points in $\E_X$ to any of the four points
in $\E_Y$, and see that each of these lengths is less than $\ell(\g|_{[t_1,t_2]})$.)
\end{remark}

\subsection{Geodesic antipodality and two types of non-uniqueness}
\label{sect4.4}

As noted in Section \ref{sect3.1}, for all $X,Y\in \sympp$ there always exists
an MSSR curve from $X$ to $Y$, the projection of some $F$-minimal geodesic.
{\em A priori}, different $F$-minimal geodesics could project to the same
MSSR curve or to different MSSR curves. It is natural to ask: Under what
conditions on $(X,Y)$ is there a unique MSSR curve from $X$ to $Y$? When
uniqueness fails, {\em how} does it fail, and what can we say
about the set $\M(X,Y)$?

For uniqueness to fail for given $X,Y$, there must be distinct $F$-minimal
geodesics $\g_i:[0,1]\to M(p)$, whose endpoints are minimal pairs
\linebreak $((U_i,D_i),(V_i,\L_i)) \in \E_X \times \E_Y ,$ $i=1,2$, such that
$F\circ\g_1\neq F\circ \g_2$. The ``how" question
above concerns the following two possibilities (not mutually exclusive):

\begin{enumerate}
\item ``Type I non-uniqueness'':  There exist such $\g_i$ whose endpoints are {\em distinct} minimal pairs $((U_i,D_i),(V_i,\L_i))$.

\item ``Type II non-uniqueness'': There exist such $\g_i$ whose endpoints are {\em the same} minimal pair $((U,D),(V,\L))$.

\end{enumerate}

\noi Since for any $D,\L\in \D^+(p)$ the minimal geodesic from $D$ to $\L$ is unique, Type II
non-uniqueness with minimal pair $((U,D),(V,\L))$ is equivalent to the existence of two or more
minimal geodesics from $U$ to $V$, which is equivalent to each of $U,V$ being in
the cut-locus (in $\Sop$) of the other.  It will be convenient for us to have some other terminology for such pairs:

\begin{defn}\label{geodantip} {\rm
Call a pair of points $(U,V)$ in $\Sop\times \Sop$ {\em geodesically antipodal}
if one point is in the cut-locus of the other (equivalently, if each point is in the
cut-locus of the other)
and {\em geodesically non-antipodal} otherwise.
Call a pair of points $((U,D), (V,\L))$ in $M(p)\times M(p)$ geodesically antipodal if
$(U,V)$  is a geodesically  antipodal pair in $\Sop\times \Sop$, and geodesically non-antipodal otherwise.
}
\end{defn}

As mentioned earlier, the cut-locus of the identity $I\in\Sop$ is precisely the set of all involutions in $\Sop$. Furthermore, because of the invariance of the Riemannian metric
$\gsop$,  an element $V\in \Sop$ is in the cut-locus of element $U$ if and
only if $V^{-1}U$ is in the cut-locus of $I$. Note that,  as would be true in
any group, if any of the elements
$V^{-1}U, UV^{-1}, U^{-1}V, VU^{-1}$ is an involution, so are all the others.

Note that a pair $(U,V)$ in $\Sop$ can be geodesically antipodal without
either point being  maximally remote from the other.  (For example, with
$p=4$, the matrix ${\rm diag}(-1,-1,1,1)$ is an involution, but is closer to
the identity $I$ than is the involution $-I$.) However, if $(U,V)$ is geodesically
antipodal, then there exists a (not necessarily unique) closed geodesic
in $\Sop$ containing $U$ and $V$, isometric to a circle  of some radius, such
that $U$ and $V$ are antipodal points of this circle in the usual sense.

Proposition \ref{distprop} is a starting-point for understanding the set
$\M(X,Y)$ for all $p$ and all $X,Y\in\sympp$: it assures us that, for
any $(U,D)\in \E_X$, every MSSR curve from $X$ to $Y$ corresponds to
some minimal pair whose first element lies in the connected component
$[(U,D)]$ of $\E_X$.  But even once we know all the minimal pairs, to
completely understand $\M(X,Y)$---or even just determine its
cardinality---we need a way to tell whether MSSR curves corresponding
to two (not necessarily distinct) minimal pairs with first point in $[(U,D)]$
are the same. (This is true whether the non-uniqueness, if any, in $\M(X,Y)$
is of Type I, Type II, or a mixture of both). Proposition \ref{prop:MSRequal}
below provides such a tool. This proposition, like Proposition \ref{distprop},
plays a crucial role in \cite{GJS2017eigen_EJS} (where it is stated without
proof), enabling an explicit computation of the sets $\M(X,Y)$ for $p=3$.

\begin{prop} \label{prop:MSRequal}
Let $X,Y \in\sympp, X\neq Y$. For $i=1,2$ assume that
$\chi_i=F\circ\g_i$ is a minimal smooth scaling-rotation curve from
$X$ to $Y$ corresponding to the minimal pair
$((UR_{U,i},D),(VR_{V,i}\, P_{g_i}^{-1}, \L_i)),$ where $R_{U,i}\in
G_D^0, R_{V,i}\in G_\L^0$, $g_i\in \tsp^+$, $\L_i=\pi_{g_i}\dotprod
\L$, and $\g_i:[0,1]\to M(p)$ is a geodesic. (We do not assume that the
two minimal pairs are distinct.) Then $\chi_1=\chi_2$ if and only if
the following two conditions hold.

\begin{itemize}

\item[(i)] Both pairs $(UR_{U_i}, VR_{V,i}
P_{g_i}^{-1})$ are geodesically non-antipodal and

\be\label{a2eqa1_genp}
R_{V,2}\,
P_{g_2}^{-1}
R_{U,2}^{-1} =R_{V,1}\,
P_{g_1}^{-1}
R_{U,1}^{-1}\,,
\ee

\noi or both pairs  are geodesically antipodal and

\be\label{a2eqa1_genp-b}
({\rm proj}_{\Sop}\g_1'(0))R_{U_1}^{-1} = ({\rm proj}_{\Sop}\g_2'(0))R_{U_2}^{-1},
\ee

\noi  where for any $(U',D')\in M(p)$, ${\rm proj}_{\Sop}$ denotes the
natural projection $T_{(U',D')}(M(p))\to T_{U'}(\Sop)$.

\item[(ii)] There exist $g\in\tsp^+, R\in G_{D,\L_1}^0$ such that

\begin{eqnarray}
\label{d2pdd1}
D &=& \pi_g\dotprod D, \\
\L_2 &=& \pi_g\dotprod \L_1, \\
\mbox{\rm and} \ \ \ R_{U,1}^{-1}R_{U,2} &=& R
P_g^{-1}
\,.
\label{rppp}
\end{eqnarray}

\end{itemize}

\noi Equation \eqref{a2eqa1_genp-b}
implies equation \eqref{a2eqa1_genp}, so \eqref{a2eqa1_genp} is always a necessary
condition for the equality $\chi_1=\chi_2$.
\end{prop}

In Proposition \ref{prop:MSRequal}, in the geodesically non-antipodal case we use
{\em endpoint data} to tell
whether the projections to $\sympp$ of two minimal geodesics from
$\E_X$ to $\E_Y$ are equal.  We will deduce this proposition
from the following theorem, proven in
\cite{JSG2015scarot}, that gives a criterion based on {\em initial-value data} to tell
whether the projections of two geodesics emanating from $\E_X$ are
equal.  In this theorem, $G_{D,L}:=G_D\intersect G_L$,
$\frakg_{D,L}=:\frakg_D\intersect \frakg_L$ (the Lie algebra of
$G_{D,L}$), and for $A\in \Sop$, $\ad_A:\sop\to\sop$ is the linear map
defined by $\ad_A(B)=[A,B]$.

\begin{notation}\rm  For $(U,D)\in M(p)$, $A\in \sop$, $L\in \D(p)$,
and any interval $I$ containing $0$, we write
$\g_{U,D,A,L}$ for the geodesic $I\to M(p)$ defined by $t\mapsto
(e^{tA}U, e^{tL}D)$ .
\end{notation}

\begin{thm}[{\cite[Theorem 3.8]{JSG2015scarot}}]\label{jsg1thm3.8}
For $i=1,2$ let $(U_i,D_i)\in M(p)$, $A_i\in \sop$, $L_i\in\D(p)$, and let
$\check{A}_i = U_1^{-1}A_iU_1$.
Let $I$ be a positive-length interval containing $0$.  Then the smooth
scaling-rotation curves $\chi_i:=F\circ \g_{U_i,D_i,A_i,L_i}: I\to \sympp$
are identical if and only if (i) $\check{A}_2 - \check{A}_1 \in
\frakg_{D_1,L_1}$, (ii) $(\ad_{\check{A}_2})^j(\check{A}_1) \in
\frakg_{D_1,L_1}$ for all $j \ge 1$, and (iii) there exist $R\in
G_{D_1,L_1}$ and $g\in \tsp^+$, such that $U_2=U_1RP_g^{-1}$,
$D_2=\pi_g\dotprod D_1$, and $L_2=\pi_g\dotprod L_1$.\footnote{In
  \cite[Theorem 3.8]{JSG2015scarot}, $g$ was actually required to be
 a particular pre-image of $\pi$ in $\tsp^+$,
but the same argument as in the proof of Proposition \ref{jsg1thm3.3}
of the present paper shows that this restriction can be removed.}
\end{thm}

To deduce Proposition \ref{prop:MSRequal} from Theorem \ref{jsg1thm3.8},
we first prove two lemmas. Beyond helping us to prove the Proposition,
these lemmas may be useful in future analysis of MSSR curves. In these
lemmas, for any $X\in\sympp$ we write $\frakg_X$ for the Lie algebra
of the stabilizer $G_X:=\{ U\in G : UXU^T = X\}$; thus $\frakg_X=\{A\in\so(p) :
AX=XA\}$. (Observe that the notation $G_X$ is consistent with the notation
$G_D$ introduced earlier for diagonal matrices.)

\begin{lemma}\label{perpatendpts}
Let $X,Y\in\sympp$ and suppose that $\chi:[0,1]\to \sympp$ is a
minimal smooth rotation-scaling curve with $X:=\chi(0)\neq Y:=\chi(1)$. Let
$\g=\g_{U,D,A,L}:[0,1]\to \Sop\times{\rm
  Diag}^+(p)$ be a geodesic for which $\chi=F\circ\g$. Then $A\in
(\frakg_X)^\perp\intersect(\frakg_Y)^\perp$, where the orthogonal
complements are taken in $\so(p)$.
\end{lemma}

\pf Since $\g$ is a smooth curve of minimal length connecting the
submanifolds $\E_X$ and $\E_Y$ of $M(p)$, the velocity vectors
$\g'(0), \g'(1)$ must be perpendicular to the tangent
spaces $T_{\g(0)}\E_X, T_{\g(1)} \E_Y$, respectively
(\cite[Proposition 1.5]{CE1975}). Making natural tangent-space
identifications, we have $T_{\g(0)}\E_X=T_{(U,D)}\E_X =
U\frakg_D\plus\{0\} \subset U\frakg_D\plus\D(p)$, where
$U\frakg_D:=\{UC : C\in \frakg_D\}$. Let $\check{A}=U^{-1}AU$. Since
$\g'(0)=(U\check{A}, DL)$, and the Riemannian metric we are using on
$\So(p)$ is left-invariant, the condition $\g'(0)\perp T_{\g(0)}\E_X$
is equivalent to $\check{A}\in(\frakg_D)^\perp$, hence to $A\in
U(\frakg_D)^\perp U^{-1}$.  Using additionally the right-invariance of
the metric on $\so(p)$, we have $U(\frakg_D)^\perp U^{-1} = (U\frakg_D
U^{-1})^\perp$. From general group-action properties, it is easily
seen that $U\frakg_D U^{-1} = \frakg_{UDU^{-1}}$. Since $UDU^{-1}=X$,
it follows that $A\in (\frakg_X)^\perp$. A similar argument at the
point $(V,\L):=\g(1)$ shows that $A\in (\frakg_{V\L V^{-1}})^\perp =
(\frakg_Y)^\perp$. \qedns

\begin{lemma}\label{MSRequal} In
  the setting of Theorem \ref{jsg1thm3.8}, assume that the smooth
  scaling-rotation curve $\chi_1$ is minimal. Then conditions (i) and
  (ii) in the theorem can be replaced by the single condition $A_2=A_1$.
\end{lemma}

\pf With notation as in Theorem \ref{jsg1thm3.8}, assume that
$\chi_2=\chi_1$.  Then the Theorem implies that $U^{-1}(A_2-A_1)U\in
\frakg_{D,L}\subset \frakg_D,$ implying that $A_2-A_1\in U\frakg_D
U^{-1} = \frakg_X$ (as in the proof of Lemma \ref{perpatendpts}). But
since $\chi_1$ is minimal, Lemma \ref{perpatendpts} implies that both
$A_2$ and $A_1$ lie in $(\frakg_X)^\perp$, hence that
$A_2-A_1\in(\frakg_X)^\perp$. Hence $A_2-A_1=0$, i.e. $A_2=A_1$.

Conversely, assume that $A_2=A_1$. Then conditions (i) and (ii) are
satisfied trivially. \qedns

\noi {\bf Proof of Proposition \ref{prop:MSRequal}:} For $i\in\{1,2\}$
let $U_i= UR_{U,i},\ V_i=VR_{V,i} P_{g_i}^{-1},$ and $\L_i=
\pi_{g_i}\dotprod \L_i\,.$

By hypothesis $\chi_i=F\circ \g_i$, where $\g_i=\g_{U_i,D,A_i,L_i}:[0,1]\to M(p)$
(for some $A_i\in\sop,L_i\in\D(p)$) is a minimal geodesic from $(U_i,D)$ to
$(V_i,\L_i)$.  Hence $L_i=\log(\L_i D^{-1})$ and $A_i\in\log(V_iU_i^{-1})$ (we write ``$\in$'' rather than
``$=$'' since if  $R$ is an involution, `` $\log R$'', as we have defined it, is a set with more
than one element; see Section \ref{sect3.1}).

It is straightforward to show that $G_{D,L_i}=G_{D,\L_i}$.  From Lemma
\ref{MSRequal}, the conditions (i) and (ii) in Theorem
\ref{jsg1thm3.8} in the equality-conditions for $\chi_1$ and $\chi_2$
can be replaced by the single condition $A_2=A_1$.

If $A_2=A_1$ then $V_2U_2^{-1}=V_1U_1^{-1}$, implying that either both pairs
$(U_i,V_i)$
are geodesically antipodal or both are geodesically non-antipodal.
In the converse direction, suppose that the pairs $(U_i,V_i)$ are geodesically non-antipodal and that $V_2U_2^{-1}=V_1U_1^{-1}$. Then  $A_2=\log(V_2U_2^{-1}) = \log(V_1U_1^{-1}) =A_1$.
Whether or not the pairs $(U_i,V_i)$ are geodesically antipodal, by definition
$({\rm proj}_{\Sop}\g_i'(0))U_i^{-1}=A_i$,  so if \eqref{a2eqa1_genp-b} holds then $A_2=A_1$.
Hence the condition $A_2=A_1$ is equivalent to condition (i) in Proposition \ref{prop:MSRequal}.

Next, letting $D$ play the role of $D_1$ in Theorem
\ref{jsg1thm3.8}, condition (iii) in the Theorem is equivalent to the existence of $g\in\tsp^+,
R\in G_{D,\L_1}$ such that $D =\pi_g\dotprod D$, $L_2 =\pi_g\dotprod
L_1$, and $U_2 =U_1 R P_g^{-1}$. But for all such $R, \pi$, we have
$RP_g^{-1}=R_0P_{g_0}^{-1}$ for some $R_0\in G_{D,\L_1}^0$ and
$g_0\in\tsp^+$ with $\pi_{g_0}=\pi_g$. Furthermore, for any $\pi\in
S_p$, if $\pi\dotprod D=D$ then $L_2=\pi\dotprod L_1 \iff \L_2
=\pi\dotprod \L_1$.  Hence, under the hypotheses of  Proposition \ref{prop:MSRequal}, condition (iii)
in Theorem
\ref{jsg1thm3.8} is equivalent to condition (ii) stated in the Proposition.

This establishes the ``if and only if'' statement in the Proposition.  The final statement of the proposition
follows from the fact that, in the notation of this proof, \eqref{a2eqa1_genp-b} is the equality $A_2=A_1$
(after multiplying both sides of \eqref{a2eqa1_genp-b} on the right by $U^{-1}$),
an equality that implies $V_2U_2^{-1}=\exp(A_2)=\exp(A_1) = V_1U_1^{-1} $.
\qedns

\subsection{Type I and Type II non-uniqueness}\label{sect:uniqq}

Within the scaling-rotation framework, the motivation to understand Type II
non-uniqueness is its effect on a true scaling-rotation {\em metric} $\rho_\SR$
on $\sympp$, mentioned earlier, that we construct from $\dsr$ in
\cite{GJS2016metric}. Various constructions and assertions concerning this
metric are simplified when we know that Type II non-uniqueness does not occur.
But, as we shall see, the study of Type II non-uniqueness also leads to geometric
results outside the scaling-rotation framework.

For small enough values of $p$, Type II non-uniqueness never occurs; for large enough $p$, it always occurs  (see Corollaries \ref{pleq4cor} and \ref{gamp} below).
Our main tool for ruling out Type II non-uniqueness
is based on a property we call {\em sign-change reducibility} (for want of a better term), defined shortly.

To motivate the definition, let $X,Y\in \sympp$ and let $((U,D),(V,\L))\in \E_X\times \E_Y$ be a minimal pair. Then one minimizer $(g,R_U,R_V)$ of the expression in brackets on the right-hand side of \eqref{fibdist8a}
is the triple $(e,I,I)$, where $e$ is the identity element of $\tsp^+$. Hence
for all $g\in \tsp$ with $\pi_g\dotprod\L = \L$---i.e. for all $g\in
\tilde{K}_{\J_\L}$ (see Notation \ref{defkj})---we must have  $\dsop(UP_g,V)=\dsop(U,VP_g^{-1})\geq \dsop(U,V).$
But $\I_p^+\subset \tilde{K}_\J$ for all $\J$, so, in particular, we must have  $\dsop(UI_\mbs,V)\geq \dsop(U,V)$ for all $\boldsig\in\I_p^+$.

\begin{defn}\label{defscr}{\rm Call a pair of points $(U,V)\in \Sop\times \Sop$
{\em sign-change reducible} if $\dsop(UI_\mbs, V)<\dsop(U,V)$ for some $\boldsig\in \I_p^+$.}
\end{defn}

From the discussion preceding Definition \ref{defscr}, we have the following:

\begin{cor} Let $((U,D),(V,\L))\in M(p)\times M(p)$. If $(U,V)\in \Sop\times\Sop$ is sign-change reducible, then $((U,D),(V,\L))$ is not a minimal pair. \qedns
\end{cor}

Sign-change reducibility is studied in more detail in Sections
\ref{sect:app}--\ref{sec:app.5}; a long digression from the topic
of scaling-rotation distance and MSSR curves is needed (but has bonuses).
Below, we summarize some results proven there, and their consequences.
Two of the main results are given in
the following Proposition (proven in Section \ref{sec:app.5}):

\begin{prop}\label{pleq4prop} (a) For $p\leq 4$, every geodesically antipodal pair $(U,V)$ in
\linebreak $\Sop\times \Sop$ is sign-change reducible. (b) For $p\geq 11$, there exist geodesically antipodal pairs $(U,V)$ in $\Sop\times \Sop$ that are {\em not} sign-change reducible.
\end{prop}

Thus the largest dimension $p_1$ for which every geodesically antipodal pair $(U,V)$ in
$\So(p_1)\times \So(p_1)$ is sign-change reducible satisfies $4\leq p_1\leq 10$.  A combination of theory and numerical
evidence leads the authors to believe that $p_1$ is closer to 10 than to 4.

An immediate consequence of Proposition \ref{pleq4prop} (a) is the following. (Again, we do not
believe the number ``4'' here
is sharp.)

\begin{cor}\label{pleq4cor}
 For $p\leq 4$, every minimal pair in $M(p)\times M(p)$ is geodesically non-antipodal.
Hence for $p\leq 4$,  for all
$X,Y\in \sympp$ for which $|\M(X,Y)|>1$, the non-uniqueness is purely of Type I.
\end{cor}

 Part of the importance of sign-change reducibility comes from the following:

\begin{prop}\label{nscrminimal} Suppose that $(U,V)$ is a pair in $\Sop\times \Sop$ that is
 not sign-change reducible.  Then there exist $D,\L\in {\cald}_{\J_{\rm
   top}}$ such that the pair $((U,D), (V,\L))$ is minimal.
\end{prop}

We will prove this below. But first note that an immediate corollary of Propositions \ref{pleq4prop}(b)  and \ref{nscrminimal} is:

\begin{cor}\label{gamp} For $p\geq 11$, there exist geodesically antipodal, minimal pairs
$((U,D),(V,\L)) \in \S_{\J_{\rm top}}\times\S_{\J_{\rm top}}\subset M(p)\times M(p)$.  Hence, for $p\geq 11$, there exist $X,Y\in \S_{[\J_{\rm top}]}\subset \sympp$ for which the set
$\M(X,Y)$ exhibits Type II non-uniqueness.
\end{cor}

Thus sign-change reducibility is more than an {\em ad hoc} criterion for
ruling out Type II
non-uniqueness for small enough $p$.  Proposition \ref{nscrminimal} and Corollary \ref{gamp}
show that, in some sense, sign-change reducibility is the {\em only} obstruction to having points
$X,Y$ in the top stratum of $\sympp$ for which $\M(X,Y)$ exhibits Type II non-uniqueness.

For $X$ or $Y$ not in the top stratum of $\sympp$, the relationship between Type II non-uniqueness and sign-change reducibility of minimal pairs in $\E_X\times \E_Y$ situation is more complicated to
analyze.  We do not investigate this relationship further in this paper.

To prove Proposition \ref{nscrminimal} we start with a lemma:

\begin{lemma}\label{gaplemma} Let $c>0$. There exist $D,\L\in \cald_{\rm top}:=\cald_{\J_{\rm
   top}}$ such that \linebreak $\|\log (D^{-1}(\pi\dotprod\L)\|^2 >
  \|\log (D^{-1}\L)\|^2+c$ for all non-identity $\pi\in S_p$.
\end{lemma}

\pf Let $c_1=\sqrt{c/(3p)}$ and let $\{a_i\}_{i=1}^p$ be a sequence of
numbers satisfying $a_{i+1}-a_i>(2\sqrt{p}+1)c_1$ for $1\leq i\leq
p-1$.  Then $|c+a_j-a_i|>2\sqrt{p}c$ for all $i\neq j$.  Let $D={\rm
  diag}(e^{a_1},\dots, e^{a_p})$ and let $\L=e^{c_1} D$.  Then
$D,\L\in \cald_{\rm top}$ and $\|\log (D^{-1}\L)\|^2 =\|c_1 I\|^2 =
pc_1^2.$

Let $\pi\in S_p, \pi\neq {\rm id},$ and let $i$ be such that
$\pi^{-1}(i)\neq i$.  Then

\ben
\|\log (D^{-1}(\pi\dotprod\L)\|^2 \geq
|c_1+a_{\pi^{-1}(i)}-a_i|^2> (2\sqrt{p}\, c_1)^2 = \|\log (D^{-1}\L)\|^2
+c.
\een
\qedns

 \ssn {\bf Proof of Proposition \ref{nscrminimal}}.  Let $D,\L\in \cald_{\rm top}$ be such that

\be\label{gap}
\|\log (D^{-1}(\pi\dotprod\L)\|^2 > \|\log (D^{-1}\L)\|^2 +
\linebreak k\,{\rm diam}(\Sop)^2
\ee

\noi for
all non-identity $\pi\in S_p$; such $D,\L$ exist by Lemma \ref{gaplemma}.
Let $X=F(U,D), Y=F(V,\L)$. The subgroups $G_D^0, G_{\L}^0$ of
$\Sop$ are trivial,
as are the subgroups
$\Gamma_{\J_D}^0$ and $\Gamma_{\J_{\L}}^0$ of $\tsp^+$. Hence in
Proposition \ref{distprop} we have $Z=\tsp^+$ and
\bearray\nonumber
\dsr(X,Y)^2 &=&
\min_{g\in \tsp^+}
  \left\{
    k\,\dsop\left(U, VP_g^{-1}\right)^2
   +
   \|  \log\left(D^{-1}(\pi_g\dotprod \L)\right)\|^2
      \right\}\\
      \nonumber
&=&
\min_{\pi\in S_p}
\left\{
k\,\min\left\{
 \dsop\left(U, VP_g^{-1}\right)^2 : g\in \tsp^+, \pi_g=\pi\right\}\right.
\\
\nonumber
&&
\left.\phantom{
\left\{\dsop\left(U, VP_g^{-1}\right)^2 \right.}
\mbox{\hspace{-.5in}}
 +
\|  \log\left(D^{-1}(\pi\dotprod \L)\right)\|^2
 \right\}.
 \nonumber
\label{distcompar2}
\eearray
\noi For all non-identity $\pi\in S_p$ and all $g_1,g_2\in\tsp^+$ with
$\pi_{g_1}=\id$ and $\pi_{g_2}=\pi$, using \eqref{gap} we then have
\bestar
d_M((U,D), g_1\dotprod (V,\L))^2
&=& k\,\dsop\left(U, V P_{g_1}^{-1}\right)^2
 +\|  \log\left(D^{-1}\L\right)\|^2
 \\
 &\leq &k\,{\rm diam}(\Sop)^2+\|  \log\left(D^{-1}\L\right)\|^2 \\
 &<& \|  \log\left(D^{-1}\left(\pi\dotprod \L\right)\right)\|^2 \\
 &\leq& d_M((U,D),
g_2\dotprod (V,\L))^2.
 \eestar
 \noi Hence the identity permutaton is the only element of $S_p$ for
 which the expression inside the outer braces in \eqref{distcompar2}
 achieves the minimum over all $\pi\in S_p$. But $\{g\in\tsp^+
 :\pi_g=\id\}$ is precisely the sign-change subgroup $\I_p^+$, and by
 hypothesis $(U,V)$ is not sign-change reducible. Hence
\bestar
\dsr(X,Y)^2 &=&
\min_{\mbs\in \I_p^+}
  \left\{
    k\,\dsop\left(U, V
I_\mbs
\right)^2
   +
   \|  \log\left(D^{-1}\L\right)\|^2
      \right\}
      \\
&=&
 k\,\dsop(U, V)^2
   +
   \|  \log\left(D^{-1}\L\right)\|^2 \\
&=&d_M((U,D),(V,\L))^2.
\eestar
\noi Thus $((U,D),(V,\L))$ is a minimal pair. \qedns

\setcounter{equation}{0}
\section{Involutions, sign-change reducibility, and distance between subspaces of $\bfr^p$ }\label{sect:app}

In this section we begin our study of sign-change reducibility.  This culminates in Section \ref{sec:app.5} with
the proof of Proposition \ref{pleq4prop}  (which, as we have seen, implies Corollary \ref{gamp}, our main result concerning Type II
non-uniqueness), but
we discover some other interesting facts along the way.
As we shall see,
questions concerning the seemingly {\em ad hoc} notion of sign-change reducibility can be translated
into questions about distances between subspaces of $\bfr^p$; for example, Proposition \ref{equivstat} states the equivalence between a sign-change-reducibility question and a question purely about the geometry of the Grassmannian $\Gr_m(\bfr^p)$ (endowed with a standard metric).
Thus, some unexpected benefits
of our investigation of Type II non-uniqueness are results, possibly of independent interest,
concerning the geometry of Grassmannians and, more generally, principal angles between subspaces of
$\bfr^p$.

Since
$\dsop(U,V)=\dsop(V^{-1}U,I)$ for $U,V\in\Sop$, the set of distances between geodesically antipodal
points in $\Sop$ is the same as the set of distances between the
identity and involutions.  Thus to understand which (if any) geodesically antipodal pairs $(U,V)$ in
$\Sop$ are sign-change reducible, it suffices to study the case $(U,V)=(R,I)$, where $R$ is an
involution.

\begin{defn} {\rm
$~$

\begin{enumerate}

\item Call $R\in \Sop$ {\em sign-change reducible} if $\dsop(RI_\mbs, I)<\dsop(R,I)$ for some
$\boldsig\in \I_p^+$ (equivalently, if the pair $(R,I)$ is sign-change reducible).  Note that sign-change reducibility of the pair $(U,V)$, as previously defined in  Definition \ref{defscr}, is equivalent to sign-change reducibility of $V^{-1}U$.

\item For $\boldsig=(\sigma_1,\dots, \sigma_p)\in\I_p$, define the {\em level} of
$\boldsig$, written $\lev(\boldsig)$, to be $\#\{i :\sigma_i=-1\}$.

\item For any involution $R\in\Sop$,
define  the {\em level} of $R$, written $\lev(R)$,
to be $\dim(E_{-1}(R)),$ where $E_{-1}(R)$ is the $(-1)$-eigenspace of $R$.
We write ${\rm Inv}(p)$ for the set of involutions in $\Sop$, and
for $0<m\leq p$ we write ${\rm Inv}_m(p)$ for the set of involutions in
$\Sop$ of level $m$. Note that $\dim(E_{-1}(R))$ is even for any $R\in \Sop$, so
${\rm Inv}_m(p)$ is empty unless $m$ is even and at least 2. Thus ${\rm Inv}(p)=
\Union_{{\rm even}\ m \geq 2}{\rm Inv}_m(p)$ (a disjoint union).

\item Let $R\in\Sop$ be an involution.  We say that $R$ is {\em reducible by a sign-change
of level $m$} if there exists $\boldsig\in \I_p^+$ of level $m$ such that $\dsop(RI_\mbs, I)
\linebreak <\dsop(R,I)$.

\end{enumerate}
}
\end{defn}

Observe that for non-identity $\boldsig\in\I_p^+$, the matrix $I_\mbs$ is an involution in $\Sop$, and  $\lev(\boldsig)=\lev(I_\mbs)$.

\begin{remark}[Involutions and Grassmannians]\label{grassremark}{\rm
The space ${\rm Inv}(p)$ can be naturally identified with a disjoint union of Grassmannians,
because
an involution $R\in\Sop$ is completely determined by its $(-1)$-eigenspace
$E_{-1}(R)$. Let
$\Gr_m(\bfr^p)$ denote the Grassmannian of $m$-planes in $\bfr^p$, and for even $m\in (0,p]$ define
$\Phi_{m,p}:\Gr_m(\bfr^p)\to {\rm Inv}_m(p)$ to
 be the map carrying $W\in \Gr_m(\bfr^p)$ to the involution in $\Sop$ whose $(-1)$-eigenspace
 is $W$.  (Thus $E_{-1}(R)=\Phi_{m,p}^{-1}(R)$ for all $R\in {\rm Inv}_m(p)$.)
Concretely, letting $\pi_V: \bfr^p\to V$ denote orthogonal projection onto
any subspace $V,$ and letting $P_V$ denote the matrix of $\pi_V$ with respect to
the standard basis of $\bfr^p$, the map $\Phi_{m,p}$ is given by

\be\label{defPhi}
\Phi_{m,p}(W)=P_{W^\perp}-P_W =
I-2P_W,
\ee

\noi reflection about the $(p-m)$-plane $W^\perp$. It is not hard to
show that ${\rm Inv}_m(p)$
is a submanifold of $\Sop$ and that $\Phi_{m,p}$ is a diffeomorphism from
$\Gr_m(\bfr^p)$ to this submanifold.
}
\end{remark}

Our study of sign-change reduciblity of involutions will make frequent use of the {\em normal form} of
an element of $\Sop$, so we review this before proceeding.

\subsection{Normal form and
distance to the identity in $\Sop$}\label{sec:app.1}

Let $k=\floor{\frac{p}{2}}$. Recall that every $R\in\Sop$ has a {\em
  normal form}: a block-diagonal matrix that, for $p$ even, is of the
form

\be\label{normform}
\sfr(\th_1, \dots, \th_k)=
\left[
\begin{array}{cccccc}
C(\th_1) &&&&& \\
& C(\th_2) &&&&\\
&& . &&&\\
&&& . &&\\
&&&& . &\\
&&&&& C(\th_k)
\end{array}
\right],
\ee

\noi where

\be\label{cth}
C(\th)=\left[\begin{array}{cc} \cos\th & -\sin\th \\ \sin\th & \cos\th
\end{array}\right]
\ee

\noi and where $\th_i\in [0,\pi], 1\leq i\leq k$.  (This can be derived quickly from the
normal form of an antisymmetric matrix, since the compactness of $\Sop$ guarantees that the exponential map $\sop\to\Sop$ is onto.) For the odd-$p$
case, the normal-form matrix is the matrix \eqref{normform} with one
more row and column appended, and with a 1 in the lower right-hand
corner (and zeroes everywhere else in the last row and column).  In
this case we define $\th_{k+1}=0$, so that for both even and odd $p$
we can
use the notation $\sfr(\th_1, \dots, \th_{\ceil{p/2}})$ for the normal
form.

Note that

\be\label{logcth}
C(\th) = \exp(\th J)\ \ \mbox{where}\ \
J=\left[\begin{array}{cc} 0 &
    -1 \\ 1 & 0
\end{array}\right].
\ee

\noi For each $R\in\Sop$ there exists an orthonormal basis of $\bfr^p$
with respect to which the linear transformation $\bfr^p\to \bfr^p$,
$v\mapsto Rv$, has matrix $\sfr(\th_1,\dots, \th_{\ceil{p/2}})$.
Thus
there exists $Q\in O(p)$ such that

\be\label{conjnorm}
R=Q \sfr(\th_1,\dots, \th_{\ceil{p/2}}) Q^{-1}.
\ee

The normal form of a given $R$ is unique up to ordering of the blocks;
the multi-set $\{\th_1,\dots, \th_{\ceil{p/2}}\}$ is uniquely
determined by $R$. From \eqref{cth} and \eqref{conjnorm} we have

\be \label{conjlognorm}
R=Q\exp(A(\th_1,\dots, \th_{\ceil{p/2}}))Q^{-1}
=\exp(QA(\th_1,\dots, \th_{\ceil{p/2}})Q^{-1})
\ee

\noi where $A(\th_1,\dots,
\th_{\ceil{p/2}})$ is the block-diagonal matrix obtained by replacing
$C(\th_i)$ by $\th_i J$ in \eqref{normform}, $1\leq i\leq \floor{p/2}$, and,
in the odd-$p$ case, replacing the 1 in the lower right-hand corner by
0. Since the normal form is unique up to block-ordering, it follows
that

\be\label{dri}
\dsop(R,I) ^2= \sum_{i=1}^{\floor{p/2}} \th_i^2 = \sum_{i=1}^{\ceil{p/2}} \th_i^2
\ee

\noi Furthermore, from \eqref{conjnorm} and \eqref{cth} it follows that

\be\label{rsym}
\rsym := \frac{R+R^T}{2} =
Q\left[
\begin{array}{cccccc}
\cos\th_1\ I_{2\times 2} &&&&& \\
& \cos\th_2\ I_{2\times 2}  &&&&\\
&& . &&&\\
&&& . &&\\
&&&& . &\\
&&&&& \cos\th_k\ I_{2\times 2}
\end{array}
\right]Q^{-1}
\ee

\noi if $p$ is even; for odd $p$ we again just append one more row and
column of the middle matrix, with a 1 in the lower right-hand corner.
Hence the values $\cos\th_i$ (and therefore the values $\th_i\in [0,\pi]$) can be recovered from
$R$ as the eigenvalues of $\rsym$, with the multiplicity of an
eigenvalue $\l$ of $\rsym$ equal to twice the multiplicity $m_\l$ of
$\l$ in the list $\cos\th_1,\dots, \cos\th_k$ in the even-$p$ case;
for odd $p$ the only difference is that multiplicity of the eigenvalue
1 of $\rsym$ is $2m_1+1$ .

\begin{remark}[Normal form, involutions, and distances to identity]\label{rem:invonorm}
$~$

\noi {\rm
Writing $R\in\Sop$ in the form
\eqref{conjnorm}, it is easily seen that $R$ is an involution if and
only if (i) for each $i$, $\theta_i$ is either $0$ or $\pi$, and
(ii) $\th_i=\pi$ for at least one $i$. For such $R$, if $\th_i=\pi$ for exactly $m$
values of $i$, then $\|A(\th_1,\dots, \th_{\ceil{p/2}})\|^2 =
m\pi^2$.
Hence if $R\in\Sop$ is an involution of level $m$, then
\be\label{invodist2}
\dsop(R,I)^2 = \frac{m}{2}\pi^2\ .
\ee
\noi Thus

\be\label{invodist}
\left\{\dsop(R,I): R\in \Sop, R\ \mbox{an involution}\right\} =
\left\{\sqrt{m}\pi: 1\leq m\leq \left\lfloor{\frac{p}{2}}\right\rfloor\right\}.
\ee

Using \eqref{conjlognorm} it can also be shown that for
every non-involution $R\in\Sop$, there is a unique $A\in\sop$ of
smallest norm such that $\exp(A)=R$.}
\end{remark}

\begin{notation}$~$
{\rm
\begin{itemize}
\item[1.] Given $R\in \Sop$ and angles $\th_1,\dots, \th_{\ceil{p/2}}
  \in [0,\pi]$ for which \linebreak
$\sfr(\th_1,\dots,\th_{\ceil{p/2}})$ is a normal form of $R$, we
define ``redundant normal-form angles" $\tilde{\th}_i\in [0,\pi]$, $1\leq i\leq p$, by

\be\label{deftth}
\tilde{\th}_{2i-1}=\tilde{\th}_{2i}=\th_i, \\ 1\leq i\leq k=
\left\lfloor{\frac{p}{2}}\right\rfloor; \ \ \tilde{\th}_p=0\ \
\mbox{if}\ p=2k+1.
\ee

\item[2.]
For any square matrix $A$
 we write $E_\l(A)$ for the $\l$-eigenspace of $A$.

\end{itemize}

}
\end{notation}

Note that \eqref{dri} can now be written as

\be\label{dri2}
\dsop(R,I) ^2= \frac{1}{2}\sum_{i=1}^p\tilde{\th_i}^2.
\ee

\subsection{Sign-change reducibility, distances in Grassmannians, and a half-angle relation}
\label{sec:app.2}

In this section we state and discuss several results, but defer their proofs to later sections.

For $p\leq 4$ one can show, without appealing to Proposition \ref{involprop-weak} below, that every involution in $\So(p)$ is sign-change reducible.  (This sign-change redubility holds for trivial reasons for when $p=2;$ holds for slightly less trivial reasons, mentioned later in Remark \ref{noninvol}, for $p=3$; and can be shown to be hold for $p=4$ using a quaternionic approach.)
It is reasonable to wonder
whether this holds for {\em all} $p$:
\begin{question}\label{involprop}  Let $p\geq 2$.
Is every involution in $\So(p)$ sign-change reducible?
\end{question}

 Our motivation for this question is not just generalization for its own sake, however.
Potential Type II non-uniqueness complicates several aspects
of the analysis of scaling-rotation distance and the
associated metric $\rho_{\SR}$ studied in \cite{GJS2016metric}.
To understand whether the ``Type II non-uniqueness" defined in Section \ref{sect:uniqq} can occur, we need to know whether a geodesically antipodal pair in $M(p)$
can be minimal. (As discussed in Section \ref{sect:uniqq}, a geodesically non-antipodal minimal pair in $M(p)$ uniquely determines an MSSR curve in $\sympp$.) A sufficient condition for any pair $((U,D),(V,\L))$ in $M(p)\times M(p)$ to be {\em non}-minimal is that the pair $(U,V)\in \Sop$ be sign-change reducible.
Since sign-change reducibility of involutions rules out the possibility of Type II non-uniqueness,
and all involutions are sign-change reducible
for
$p\leq 4$,  it is natural to ask Question \ref{involprop} and wish for
 the answer to be yes.

 The answer, however, is more complicated.
We shall see that the answer to Question \ref{involprop} is yes for $p\leq 4$ and no for $p\geq 11$ (we do
not know the answer for $5\leq p\leq 10$), but that for all $p$,
involutions of high enough level are sign-change reducible---morevover, by a sign-change of the same level:

\begin{prop}\label{involprop-weak}
Let $R\in\Sop$ be an involution for which $\lev(R) \geq \frac{1}{2}p$.
 Then there exists
  $\boldsig\in \I_p^+$, with $\lev(\boldsig)=\lev(R)$,
 such that $\dsop(RI_\mbs, I) \linebreak < \dsop(R,I)$.
\end{prop}

 We defer the proof to Section \ref{sec:app.4}.

Since $\lev(R)=\dim(E_{-1}(R))\geq 2$ for every involution $R$, Proposition \ref{involprop-weak}
 (once proved) immediately establishes Proposition \ref{pleq4prop}(a) and Corollary \ref{pleq4cor}: for $p\leq 4$, all involutions are
sign-change reducible, and hence all minimal pairs in $M(p)\times M(p)$  are geodesically non-antipodal.

We shall see below (Proposition \ref{equivstat}) that sign-change reducibility by a sign-change of the same level is equivalent to a statement purely about the geometry of Grassmannians.
For reasons given shortly,
it seems likely to the authors that the ``same level" condition appearing in Proposition \ref{involprop-weak}
is optimal (even without the ``$\lev(R) \geq \frac{1}{2}p$" restriction)  in the sense that
$\min_{\mbs\in \I_p^+}\{\dsop(RI_\mbs, I)\}$
is achieved by a sign-change matrix $\boldsig\in\I_p^+$ for which $\lev(\boldsig)=\lev(R)$.
If this is true, then the analysis of whether an involution $R$ is sign-change reducible simplifies; we need only consider $\boldsig\in\I_p^+$ of the same level as $R$.  This (potential)
simplication is actually of greater value to us than knowing, for a given $R\in{\rm Inv}(p)$, whether all minimizers of $\dsop(RI_\mbs, I)$ have the same level as $R$, so we state only the following weaker conjecture:

\begin{conjecture}\label{onlymcanwork}
Let $m\geq 2$ be even, and let $R\in \Sop$ be an involution of level $m$. If $R$ is sign-change
reducible, then it is reducible by a sign-change of level $m$.
\end{conjecture}

In Section \ref{sec:app.4} we will prove the following special case of this conjecture:

\begin{prop}\label{conjtruem2}  Conjecture \ref{onlymcanwork} is true for $m=2$.
\end{prop}

The reason we expect more generally that
$\min_{\mbs\in \I_p^+}\{\dsop(RI_\mbs, I)\}$
is achieved by a $\boldsig$ for which $\lev(\boldsig)=\lev(R)$
is as follows.
Every sign-change matrix $I_{\mbs_1}\in\I_p^+$ is itself an involution, and satisfies
\ben
\dsop(I_{\mbs_1} I_{\mbs_1}, I)= 0 < \min\{\dsop(I_{\mbs_1} I_{\mbs}\ ,I): \boldsig\in \I_p^+, \boldsig\neq\boldsig_1\}.
\een
\noi  Thus for $R\in\Sop$ sufficiently close to
$I_{\mbs_1}$, we have
\ben
\dsop(RI_{\mbs_1},I)< \min\{\dsop(R I_{\mbs},I): \boldsig\in \I_p^+, \boldsig\neq\boldsig_1\}.
\een
\noi The function carrying an involution in $R\in\Sop$ to $\lev(R)$ is continuous,
so for $R\in{\rm Inv}(p)$ sufficiently close to
$I_{\mbs_1}$ we also have $\lev(R)=\lev(\boldsig_1)$.  Hence for every $R\in{\rm Inv(p)}$
sufficiently close to a sign-change matrix, \linebreak
$\min_{\mbs\in \I_p^+}\{\dsop(RI_\mbs, I)\}$ is achieved by a sign-change matrix having the same
level as $R$. It seems plausible that this remains true even without the ``sufficiently close to a sign-change matrix" restriction.

As noted in Remark \ref{grassremark}, for even $m\geq 2$ the space ${\rm Inv}_m(p)$ is diffeomorphic to the Grassmannian $\Gr_m(\bfr^p)$.  This Grassmannian
carries a Riemannian metric induced by Riemannian submersion from $(\Sop,\gsop)$.  It is known that the associated squared geodesic-distance between two points $W,Z\in \Gr_m(\bfr^p)$ is, up to a constant factor, simply the sum of squares of the principal angles between the two $m$-planes $W,Z$.\footnote{This fact follows from Wong's results on geodesics in \cite{Wong1967}, and has been cited elsewhere in the literature (e.g. \cite[p. 337]{EAS1998}), though the explicit statement does not appear in \cite{Wong1967}.} Choosing the normalization in which the squared geodesic distance $\dgr(W,Z)^2$ {\em equals} the sum of squares of the principal angles (equation \eqref{wongdist} below),
we will prove the following
in Section \ref{sec:app.3}:

\begin{prop}\label{isometry}
The map $\Phi=\Phi_{m,p}:  (\Gr_m(\bfr^p), \dgr)\to ({\rm Inv}_m(p),\dsop)$ (see \mbox{\eqref{defPhi}}) is an isometry, up to a constant factor of 2:

\be\label{isometry_formula}
\dsop(\Phi(W), \Phi(V))=2\dgr(W,V)
\ee

\noi for all $W,V \in \Gr_m(\bfr^p)$.

\end{prop}

We derive Proposition \ref{isometry} from a general half-angle relation proven in Section \ref{sec:app.3}:

\begin{prop}\label{cor:halfangle} Let $R_1, R_2$ be involutions in $\Sop$.  For $i=1,2$ let $m_i=\dim(E_{-1}(R_i))$, and let $m=\min\{m_1,m_2\}$.
Let $\{\th_i\in [0,\pi]\}_{i=1}^{\ceil{p/2}}$ be angles for which
$\sfr(\th_1,\dots, \th_{\ceil{p/2}})$ is a normal form of the product $R_1R_2$, and let
$\{\tilde{\th}_i\}_{i=1}^p$ be as defined in \eqref{deftth}.  Then for some injective map $\i: \{1,2,\dots,m\}
\to \{1,2,\dots,p\}$, the
principal angles between $E_{-1}(R_1)$ and $E_{-1}(R_2)$ satisfy

\be\label{relatephithetagen}
\phi_j(E_{-1}(R_1),E_{-1}(R_2)) = \frac{\tilde{\th}_{\i(j)}}{2}, \ \ 1\leq j\leq m.
\ee

\noi For every $i\notin {\rm range}(\i)$, the angle $\tilde{\th}_i$ is either $0$ or $\pi$.

\end{prop}

In other words, as stated in the introduction:  for any two involutions $R_1, R_2\in \Sop$, each of the principal angles between $E_{-1}(R_1)$ and $E_{-1}(R_2)$ is exactly half a correspondingly indexed normal-form angle of $R_1R_2$.

Proposition \ref{isometry}
can also be proven by purely Riemannian methods, but the proof
we give, via Proposition \ref{cor:halfangle}, is independent in the sense that it does not make any use of a Riemannian metric on
$\Gr_m(\bfr^p)$; see Remark \ref{rem:isom}.

 In Section \ref{sec:app.3}, after proving Proposition \ref{isometry} we will use it to deduce the
following:

\begin{prop} \label{equivstat}
Let $m,p$ be integers with $m$ even and $0<m\leq p$.
Then the following two statements  are equivalent:

\begin{enumerate}

\item Every involution
$R\in\Sop$ of level $m$ is sign-change reducible by a sign-change of level $m$.

\item For
every $W\in \Gr_m(\bfr^p)$,
there exists a coordinate $m$-plane $\bfr^J$ (see Notation \ref{jmpnotn}) such that

\be\label{grassineq}
\dgr(W,\bfr^J)^2< \frac{m\pi^2}{8} \ .
\ee

\end{enumerate}

\end{prop}

In other words, {\em the sign-change reducibility asserted in Statement 1 of the Proposition is equivalent
to a statement purely about the geometry of Grassmannians (with the metric $\dgr$)}, namely that the  coordinate $m$-planes in $\bfr^p$ form a ``lattice"
of
${p \choose m}$ points in $\Gr_m(\bfr^p)$ such that such that every point in
$\Gr_m(\bfr^p)$ is within
distance $(m\pi^2/8)^{1/2}$ of some lattice-point.  This gives us a geometric way to tackle Question \ref{involprop}, at least for sign-change reducibility of an involution $R$ by a sign-change matrix of the same level.  However, the authors do not know a formula for
$\min_{J\in {\cal J}_m}\{d_{Gr}(W,\bfr^J)\}$ for general
$W\in \Gr_m(\bfr^p)$, or (more importantly), a formula for
$\max_{W\in \Gr_m(\bfr^p)}\left\{\min_{J\in {\cal
    J}_m}\{d_{Gr}(W,\bfr^J)\}\right\}$.

Note that Proposition \ref{involprop-weak} asserts that
statement  1 of Proposition \ref{equivstat} is true
whenever $m\geq \frac{p}{2}$.  To put into perspective the number $\frac{m}{8}\pi^2$
appearing in statement  2 of Proposition \ref{equivstat}, and better understand the relevance
of the comparison between
$m$ and $\frac{p}{2}$,
note that the
squared diameter of $\Gr_m(\bfr^p)$ is $\min\{m,p-m\}\frac{\pi^2}{4}$.
So for $m\leq \frac{p}{2}$, \eqref{grassineq} is equivalent to

\be\label{grassineq2}
d_{Gr}(W,\bfr^J)^2 < \frac{1}{2}{\rm diam}(\Gr_m(\bfr^p))^2.
\ee

\noi For $m>\frac{p}{2}$, the right-hand side of \eqref{grassineq} is a greater fraction of ${\rm diam}(\Gr_m(\bfr^p))^2$, so it is ``easier" for statement 2 of Proposition \ref{equivstat} to be true for $m>\frac{p}{2}$ than for $m<\frac{p}{2}$.

\begin{remark}\label{noninvol}{\rm It is relatively easy to show that for any involution $R$,
there exists $\boldsig$ for which $RI_\mbs$ is not an involution. For
$p=2,3$, we have $\dsop(I,R)\leq \pi$ for every $R\in\Sop$, and
$\dsop(I,R)=\pi$ for every involution $R$, so any non-involution is
closer to the identity than is any involution. Hence for these values
of $p$, Proposition \ref{equivstat} is easy to prove. However, for $p\geq 4$,
given an involution $R$ and a $\boldsig\in \I_p^+$ for which $RI_\mbs$
is not an involution, \eqref{invodist} shows that we cannot
immediately deduce that $\dsop(I,RI_\mbs)<\dsop(I,R)$.}
\end{remark}

\setcounter{equation}{0}
\section{Proofs of the half-angle relation
and results related to Grassmannians: Propositions
\ref{isometry}, \ref{cor:halfangle}, and \ref{equivstat}}\label{sec:app.3}

The half-angle relation in Proposition \ref{cor:halfangle} underlies our proofs of
 of the most of the other results stated in Section \ref{sec:app.2} (all but Proposition \ref{conjtruem2}).
When the dimensions of the eigenspaces in Proposition \ref{cor:halfangle} are equal, the
half-angle relation
leads to the elegant distance-relation \eqref{isometry_formula}.  This equidimensonal case
is actually the only one we need for the application to Type-II non-uniqueness of MSSR curves. However, the half-angle relation \eqref{relatephithetagen} holds whether or not $\dim(E_{-1}(R_1)) =\dim(E_{-1}(R_2))$.
Since this fact may be of interest outside the scope of this paper,
and is not
much harder to prove without the equal-dimensions restriction, we have stated (and will prove) the more general relation.

Section \ref{sec:app.3a} is devoted to establishing Proposition \ref{cor:halfangle}.
In Section \ref{sec:app.3b}, we apply this proposition to establish
Propositions
\ref{isometry} and \ref{equivstat}.

\subsection{The half-angle relation}\label{sec:app.3a}

We start with some notation.

\begin{notation}\label{jmpnotn} $~$

{\rm
\begin{itemize}

\item[1.] For $1\leq i\leq p$ let $\bfe_i$ denote the $i^{\rm th}$
  standard basis vector of $\bfr^p$.

\item[2.]  For $0\leq m\leq p$, let $\jmp$ denote the collection of
$m$-element subsets of $\{1,\dots, p\}$.

\begin{itemize}

\item[(a)] For $0\leq m\leq p$ and $J\in\jmp$,
define $\boldsig^J=(\sigma_1, \dots, \sigma_p)\in \I_p$ by
$\sigma_i= -1$ for $i\in J$ and $\sigma_i=1$ for $i\notin J$.
Similarly, for $\boldsig=(\sigma_1, \dots, \sigma_p)\in \I_p$, define $J^\mbs=
\{i\in\{1,\dots,p\} : \sigma_i=-1\}$.  (The maps $J\mapsto \boldsig^J$ and
$\boldsig\mapsto J^\mbs$ are inverse to each other.)

\item[(b)] If $1\leq m\leq p$ and $J=\{i_1,\dots, i_m\} \in \jmp$, with $i_1<i_2<\dots <i_m$,
let ${\sf E}_J$ denote the $p\times m$ matrix whose $k^{\rm th}$ column is
$\bfe_{i_k}$, $1\leq k\leq m$.

\item[(c)] For $0\leq m\leq p$ and $J\in \jmp$, define
$\bfr^J=\{(x^1, x^2, \dots, x^p)\in \bfr^p : x^i = 0
\ \mbox{if}\ i\notin J\}$.

\end{itemize}

The collection $\{\bfr^J : J\in\jmp\}$ is the set of ``coordinate $m$-planes'' in $\bfr^p$.

\item[3.]  For any $J\subset \{1,\dots, p\}$, let $J'$ denote the complement of $J$ in
$\{1,\dots, p\}$.

\item[4.]  For $m_1,m_2\in \{1,2,\dots,p\}$, $W\in \Gr_{m_1}(\bfr^p),$
$Z\in \Gr_{m_2}(\bfr^p)$, and $J\in{\cal J}_{m_1,p}$, writing $m=\min\{m_1,m_2\}$,

\begin{itemize}
\item[(a)] let
  $\phi_1(W,Z), \dots, \phi_m(W,Z)$, denote the principal angles
  between the $m_1$-plane $W$ and the $m_2$-plane $Z$ (see
\cite[Section 12.4.3]{Golub1989}), and

\item[(b)] let $\phi_{J,i}(W)=\phi_{i}(W,\bfr^J)$, $1\leq i\leq m_1$.

\end{itemize}

\item[5.] For $1\leq m\leq p$ define $d_{Gr}:
  \Gr_m(\bfr^p)\times \Gr_m(\bfr^p)\to \bfr$ by

\be\label{wongdist}
d_{Gr}(W,Z) = \left(\sum_{i=1}^m \phi_i(W,Z)^2\right)^{1/2}.
\ee

\noi As noted earlier, $\dgr$ is the distance-function defined by the standard
$\Sop$-invariant Riemannian metric on $\Gr_m(\bfr^p)$  (up to a constant factor).

  \end{itemize}
  }
\end{notation}

The following  long but far-reaching technical lemma, giving several detailed relations
between a general involution in $\Sop$ and its product with a sign-change matrix, is our key tool
for establishing the results stated in Section \ref{sec:app.2}.   It is best thought of as
a series of lemmas, all with the same hypotheses, that have been rolled into one long lemma in
order to avoid restating hypotheses and notational definitions.
After proving the lemma, we build on it with two corollaries,
completing the groundwork for the proofs (in later sections) of the Section \ref{sec:app.2} propositions.

\begin{lemma}\label{R-eigenlemma} Let $R\in\Sop$ be an involution,
let $\boldsig\in \I_p^+$, assume $0<m_\mbs:=\lev(\boldsig)<p$, and let
$J=J^\mbs$ (see Notation \ref{jmpnotn}). Viewing
$\bfr^p$ as $\bfr^{J'}\plus \bfr^J$, below we write every $p\times p$
matrix in the block form $\left[ \begin{array}{ll} A_1 & A_2 \\ A_3 &
    A_4 \end{array} \right]$, where $A_1$ is $(p-m_\mbs)\times (p-m_\mbs)$,
$A_2$ is $(p-m_\mbs)\times m_\mbs$, $A_3$ is $m_\mbs\times (p-m_\mbs)$, and $A_4$ is
$m_\mbs\times m_\mbs$. Then:

(i) In this block form,
\be\label{rblock}
R=\left[ \begin{array}{ll} R_1 & R_2 \\ R_2^T & R_4 \end{array}
\right],
\ee
\noi where $R_1$ is a symmetric $(p-m_\mbs)\times (p-m_\mbs)$ matrix, $R_4$ is a
symmetric $m_\mbs\times m_\mbs$ matrix, and $R_2$ is $(p-m_\mbs)\times m_\mbs$.

(ii) In the same block form,

\be\label{r1r4}
(RI_\mbs)_{\rm sym} =\frac{1}{2}(RI_\mbs + I_\mbs R) =
\left[ \begin{array}{cc} R_1 & 0 \\ 0 & -R_4 \end{array}
\right].
\ee

(iii) All eigenvalues of $R_1$ and $R_4$ lie in the interval $[-1,1]$.

(iv) For every $\l\in (-1,1)$, if $\l$ is an eigenvalue of $R_1$
  (respectively, $R_4$), then $-\l$ is an eigenvalue of $R_4$
  (resp. $R_1$) with the same multiplicity.

(v)  Let $l$ denote the number of eigenvalues of $R_1$, counted with multiplicity, lying
in the interval $(-1,1)$.  Then $l$ is also the number of eigenvalues of $R_4$, counted with multiplicity, lying
in  $(-1,1)$, and $l\leq \min\{m_\mbs,p-m_\mbs\}$.

(vi) The inclusion map $\bfr^{J'}\to \bfr^p$ defined by $v\mapsto
\left[\begin{array}{l} v\\ 0 \end{array}\right]$ restricts to
isomorphisms $E_{\pm 1}(R_1)\to E_{\pm 1}(R)\intersect \bfr^{J'}$.
Similarly the inclusion map $\bfr^{J}\to \bfr^p$ defined by $w\mapsto
\left[\begin{array}{l} 0\\ w \end{array}\right]$restricts to
isomorphisms $E_{\pm 1}(R_4)\to E_{\pm 1}(R)\intersect \bfr^J$.

(vii) Let $l_- = \dim(E_1(R)\intersect \bfr^J), \ l_+ =
\dim(E_{-1}(R)\intersect \bfr^{J'} )$.\footnote{The $\pm$ subscripts
  are chosen according to the eigenspaces of $I_\mbs$ rather than $R$:
  $\bfr^J=E_{-1}(I_\mbs)$, $\bfr^{J'}=E_1(I_\mbs)$.}  Then \linebreak
  $ \dim(E_1(R_4))=l_-$ and $\dim(E_{-1}(R_1)) =l_+$.  (Thus $l_- + l_+$ is
the multiplicity of $-1$ as an eigenvalue of $(RI_\mbs)_{\rm sym}$ in \eqref{r1r4}, hence of
$RI_\mbs$ itself, and therefore yields a lower bound on $\dsop(RI_\mbs,I).$) Furthermore,

\be\label{lplm1}
l_-\geq \lev(\boldsig) - \lev(R) \ \ \
\mbox{\rm and}\ \ \  l_+\geq \lev(R)-\lev(\boldsig),
\ee

\noi and

\be\label{lplm2}
l_- - l_+ = \lev(\boldsig) -\lev(R).
\ee

 (viii) There exist an orthonormal $R_1$-eigenbasis
$\{v_i\}_{i=1}^{p-m_\mbs}$ of $\bfr^{p-m_\mbs}$ (i.e. an orthonormal basis of
$\bfr^{p-m_\mbs}$ consisting of eigenvectors of $R_1$) and an
$R_4$-eigenbasis $\{w_i\}_{i=1}^{m_\mbs}$ of $\bfr^{m_\mbs}$. For any such bases
$\{v_i\}$ of $\bfr^{p-m_\mbs}$, $\{w_i\}$ of $\bfr^{m_\mbs}$, let $\{\l'_i\}, \{\l_i\}$ be the
corresponding
eigenvalues (i.e. $R_1v_i=\l_i'
v_i$ and $R_4w_i=\l_i w_i$), and define

\bearray\label{defbfv}
\bfv_i &=&
\left\{ \begin{array}{ll}
\left[\begin{array}{c} v_i \\
\frac{1}{1+\l_i'}R_2^T v_i \end{array}\right],
&1< i\leq p-m_\mbs, \ \l_i'\neq -1,
\\ \\
\left[\begin{array}{c} v_i \\
0 \end{array}\right],
& 1\leq i\leq p-m_\mbs\ , \l_i'=-1,
\end{array}\right.
\\
\bfw_i&=&
\left\{ \begin{array}{ll}
\left[\begin{array}{c} \frac{-1}{1-\l_i}R_2w_i \\
w_i\end{array}\right],
& 1\leq i\leq m_\mbs,\ \l_i\neq 1, \\ \\
\left[\begin{array}{c} 0 \\
w_i\end{array}\right], & 1\leq i\leq m_\mbs,\ \l_i=1.
\end{array}\right.
\label{defbfw}
\eearray

\noi Then

\be\label{basisep}
\left\{
\sqrt{\frac{1+\l_i'}{2}}\bfv_i : 1\leq i\leq p-m_\mbs,\ \l_i'\neq -1
\right\} \union \left\{ \bfw_i  : 1\leq i\leq m_\mbs,\ \l_i =1 \right\}
\ee

\noi (ordered arbitrarily) is an orthonormal basis of $E_1(R)$, and the set

\be\label{basisem}
\left\{\sqrt{\frac{1-\l_i}{2}}
\bfw_i : 1\leq i\leq m_\mbs,\ \l_i\neq 1
\right\}\union \left\{ \bfv_i : 1\leq i\leq p-m_\mbs,\ \l_i' = -1\right\}
\ee

\noi (ordered arbitrarily) is an orthonormal basis of $E_{-1}(R)$.
Note that the cardinality
of the second set in \eqref{basisep} (respectively
\eqref{basisem}) is $l_-$ (resp. $l_+$).

\end{lemma}

\pf To simplify notation in this proof, we let $m=m_\mbs$.

Since $R\in\Sop$ is an involution, $R=R^{-1}=R^T$. Hence $R$ is
symmetric, implying assertion (i), and $\bfr^p$ is the orthogonal
direct sum of $E_1(R)$ and $E_{-1}(R)$ (since the only possible
eigenvalues of an involution are $\pm 1$).

For (ii), observe that in the block-form decomposition we are using,

\ben
I_\mbs=\left[ \begin{array}{cc} I_{(p-m)\times (p-m)} & 0 \\ 0 &
    -I_{m\times m}
\end{array}
\right].
\een

\noi A simple calculation then yields \eqref{r1r4}.

Next, because $R^2=I$, we have the following relations:

\bearray
\label{invo1}
R_1^2 + R_2 R_2^T &=& I_{(p-m)\times (p-m)}\ , \\
\label{invo2}
R_1R_2 + R_2R_4 &=& 0_{(p-m)\times m}, \\
\label{invo2b}
R_2^TR_1 + R_4R_2^T &=& 0_{m\times(p-m)}, \\
R_4^2 + R_2^TR_2 &=& I_{m\times m}.
\label{invo3}
\eearray

From \eqref{invo1} and \eqref{invo3}, for any $v\in
\bfr^{p-m}, w\in \bfr^m$, we have

\bearray
\label{normsum1}
\|R_1 v\|^2 + \|R_2^T v\|^2 &=& \|v\|^2,\\
\|R_4 w\|^2 + \|R_2 w\|^2 &=& \|w\|^2.
\label{normsum2}
\eearray

\noi It follows from \eqref{normsum1}--\eqref{normsum2} that if $\l$
is an eigenvalue of $R_1$ or $R_4$, then $|\l|\leq 1$, yielding (iii).

To obtain (iv), consider the operators $L:\bfr^m\to \bfr^{p-m}$ and
$L^*: \bfr^{p-m}\to\bfr^m$ defined by $L(w)=R_2 w$ and
$L^*(v)=R_2^Tv$.
Suppose that $R_1$ has an eigenvalue $\l$ with $|\l|<1$, and let
$0\neq v\in E_\l(R_1)$. Let $w=R_2^Tv$; note that \eqref{normsum1}
implies $w\neq 0$.  Using \eqref{invo2b},

\ben
R_4w = R_4R_2^Tv = -R_2^TR_1v = -R_2^T\l v = -\l w.
\een

\noi Hence $L^*$ maps $E_\l(R_1)$ injectively to $E_{-\l}(R_4)$.
Similarly, if $R_4$ has an eigenvalue $-\l$ with $|\l|<1$, and
$L^*$ maps $E_{-\l}(R_4)$ injectively to
$E_{\l}(R_1)$.

It follows that, for any $\l\in\bfr$ with $|\l|<1$, $\l$ is an
eigenvalue of $R_1$ if and only if $-\l$ is an eigenvalue of $R_4$,
and that the maps

\bearray\label{lsiso}
L^*|_{E_\l(R_1)} : E_\l(R_1)\to E_{-\l}(R_4) = E_\l(-R_4), \\
\ \ \
L|_{E_\l(-R_4)} : E_\l(-R_4) = E_{-\l}(R_4)\to E_{\l}(R_1)
\label{liso}
\eearray

\noi are isomorphisms. This establishes (iv). Statement (v) is an
immediate corollary of (iv).

For (vi), let $\i:\bfr^{J'}\to\bfr^p$ be the first inclusion map in the
lemma.  Note that $R\left[\begin{array}{l} v\\ 0 \end{array}\right] =
\left[\begin{array}{l} R_1v\\ R_2^Tv \end{array}\right]$.  If $v\in
E_\l(R_1)$ with $\l=\pm 1$, equation \eqref{normsum1} implies that
$R_2^Tv=0$, hence that $R\i(v) = \l\i(v)$.  Conversely, if $R\i(v) =
\l\i(v)$, then $R_1v=\l v$ (and $R_2^Tv=0$).  Hence $\i$ carries
$E_\l(R_1)$ isomorphically to $E_\l(R)\intersect \bfr^{J'}$. The
argument for the inclusion map $\bfr^{J}\to\bfr^p$ is essentially
identical. This establishes (vi).

Part (vi) implies that $\dim(E_1(R_4))=\dim(E_1(R)\intersect \bfr^J)= l_-$
and that $\dim(E_{-1}(R_1))=\dim(E_{-1}(R)\intersect \bfr^{J'})= l_+$,
the first assertion in (vii). To obtain \eqref{lplm1}--\eqref{lplm2},
note that for any subspaces $V, W$ of $\bfr^p$, we have

\be\label{genlinalg2}
\dim(V^\perp\intersect W) - \dim(V\intersect W^\perp)
=\dim(W)-\dim(V).
\ee

(The proof of \eqref{genlinalg2} is straightforward linear algebra.) Applying this to the case $V=E_{-1}(R), V^\perp = E_1(R),
W=E_{-1}(I_\mbs)=\bfr^J, W^\perp = E_1(I_\mbs)=\bfr^{J'}$, we have
$l_-=\dim(V^\perp\intersect W)$ and $l_+=\dim(V\intersect W^\perp)$,
so \eqref{lplm2} follows from \eqref{genlinalg2}.  The inequalities in
\eqref{lplm1} follow directly from \eqref{lplm2}.

(viii) Since $R_1$ (respectively $R_4$) is symmetric, an orthonormal
$R_1$-eigenbasis $\{v_i\}$ of $\bfr^{p-m}$ (resp., orthonormal
$R_4$-eigenbasis $\{w_i\}$ of $\bfr^m$) exists.  Select such
eigenbases, and let $\{\l_i\}$, $\{\l_i'\}$ be eigenvalues as defined
in the Lemma.  Note that the second set in
\eqref{basisep} is a basis of $E_1(R_4)$, which by (vi) is
isomorphic to $E_1(R)\intersect \bfr^{J}$. Hence the cardinality
of this set is $\dim(E_1(R)\intersect \bfr^{J})$, i.e. $l_-$.
Similarly, the second set in \eqref{basisem} is a basis of $E_{-1}(R_1)$
and has cardinality $l_+$.

Without loss of generality, we may assume that the eigenvectors $v_i$
with eigenvalue $-1$, if any, are the last $l_+$, and that the eigenvectors
$w_i$ with eigenvalue 1, if any,  are the last $l_-$.
Using \eqref{invo3}, for $1\leq i\leq m$ we have $R_2^TR_2 w_i=
(1-\l_i^2)w_i$, while
using \eqref{invo2} we find $R_1R_2
w_i
=-\l_i R_2 w_i$.
Then, using \eqref{rblock}, a
    simple calculation shows that $R\bfw_i=-\bfw_i$. Hence $\bfw_i\in
  E_{-1}(R)$ for $1\leq i\leq m-l_-$, while from part (vi), $\bfv_i\in E_{-1}(R)$
for $p-m-l_+<i\leq p-m$.

Let $\lb
\cdot\,, \cdot\rb$ denote  the standard inner product on $\bfr^n$ for any $n$.
As seen in the proof of part (vi),  $v\in E_{-1}(R_1)$ implies $R_2^Tv=0$.
Hence for $p-m-l_+<i\leq p-m$ and $1\leq j\leq m-l_-$,
$\lb \bfv_i,\bfw_j\rb \propto \lb v_i,R_2w_j\rb = \lb R_2^Tv_i,w_j\rb =0,$
while for $p-m-l_+<i,j\leq p-m$ we have $\lb \bfv_i,\bfv_j\rb = \lb v_i,v_j\rb
=\d_{ij}$. Finally, for $i,j\leq m-l_-$, using the fact that $\lb R_2 w_i, R_2 w_j\rb =
\lb w_i, R_2^TR_2 w_j\rb=\linebreak \lb w_i, (1-\l_i^2) w_j\rb$, a simple computation yields
$\lb \bfw_i, \bfw_j\rb =\frac{2}{1-\l_i}\d_{ij}\ .$
Thus \linebreak $\{\sqrt{\frac{1-\l_i}{2}} \bfw_i : 1\leq i \leq m-l_-\}
\union \{\bfv_i : p-m-l_+<i\leq m\}$ is an orthonormal
subset of $E_{-1}(R)$.  Using \eqref {lplm2}, the cardinality of this subset is
$m- l_- +l_+ =\lev(\boldsig)-(\lev(\boldsig)-\lev(R))=\lev(R) = \dim(E_{-1}(R))$.
Hence \eqref{basisem} is an orthonormal basis of $E_{-1}(R)$.

The proof that \eqref{basisep} is an orthonormal basis of $E_1(R)$ is similar.
\qedns

\begin{cor}\label{R-eigencor2}
Hypotheses and notation as in Lemma \ref{R-eigenlemma}.
Let $l_+= \linebreak
\dim(E_{-1}(R_1))$ and $l_-=\dim(E_1(R_4))$ (as in Lemma \ref{R-eigenlemma}(vii)). In addition
let $\{\th_i\in [0,\pi]\}_{i=1}^{\ceil{p/2}}$ be angles for which
$\sfr(\th_1,\dots, \th_{\ceil{p/2}})$ is a normal form of $RI_\mbs$, and let
$\{\tilde{\th}_i\}_{i=1}^p$ be as defined in \eqref{deftth}.  Let
$J_*=\{j\in J : 0<\tilde{\th}_{j}<\pi\}$.
Then $|J_*|\leq \min\{m_\mbs,p-m_\mbs\}$, and

\be\label{dri3}
\dsop(RI_\mbs,I)^2 =\frac{1}{2}(l_+ + l_-)\pi^2
+\sum_{j\in J_*} \tilde{\th}_j^2.
\ee

\noi If $\lev(\boldsig)=\lev(R)$, then

\be\label{dri5}
\dsop(RI_\mbs,I)^2 = l_-\pi^2
+\sum_{j\in J_*} \tilde{\th}_j^2 = \sum_{j\in J}\tilde{\th}_j^2.
\ee

\end{cor}

\pf
Let $\b':J'\to \{1,\dots, p-m_\mbs\}$, $\b:J\to \{1,\dots, m_\mbs\},$ be order-preserving bijections. By \eqref{rsym}, the eigenvalues of
$(RI_\mbs)_{\rm sym}$, counted with multiplicity, are
$\{\cos\tilde{\th}_i\}_{i=1}^p$.  But from Lemma
\ref{R-eigenlemma}(ii), we can read off the eigenvalues of
$(RI_\mbs)_{\rm sym}$ from \eqref{r1r4}; they are $\l_1',\dots,
\l_{p-m_\mbs}', -\l_1, \dots, -\l_{m_\mbs}$ (ordered arbitrarily).  Thus, reordering the $\l'_j$ and
the $\l_j$ appropriately,  for $1\leq j\leq p$ we have

\be\label{relatethlam}
\cos\tilde{\th}_j =\left\{ \begin{array}{ll} \l'_{\b'(j)} & \mbox{if}\ j\in J',
\phantom{\l'_{i'_{j_X}}}
\\
-\l_{\b(j)} & \mbox{if}\ j\in J.
\end{array}\right.
\ee

Define $J_*=\{j\in J: \l_{\b(j)}\neq \pm 1\}.$ Observe that $J_*$ can  also be characterized as $\{j\in J: \l_{\b(j)}\neq \pm 1\}$.
Similarly, define $J'_*=\{j\in J': \l'_{\b'(j)}\neq \pm 1\}=
\{j\in J' : 0<\tilde{\th}_{j}<\pi\}.$ By part (v) of Lemma \ref{R-eigenlemma},
$|J'_*| = |J_*|=l\leq \min\{m_\mbs,p-m_\mbs\}$, and by part (iv) of the Lemma there is a bijection
$b:J_*\to J'_*$ such
that $-\l_{j}=\l'_{b(j)}$ for all $j\in J'_*.$  Hence

\be\label{relatethlam2}
\tilde{\th}_j =\left\{ \begin{array}{ll} \cos^{-1}\l'_{\b'(j)} & \mbox{if}\ j\in J'_*, \\
\cos^{-1}\l'_{b(\b(j))} & \mbox{if}\ j\in J_*, \\
0\ \mbox{or}\ \pi & \mbox{otherwise}.
\end{array}\right.
\ee

\noi In particular,

\be\label{samesum}
\sum_{j\in J_*}\tilde{\th}_j^2 =\sum_{j\in J'_*}\tilde{\th}_j^2.
\ee

Next, note that

\be\label{lpsum}
\sum_{j\in J'\setminus J'_*} \tilde{\th}_j^2 =\sum_{\{j\in J': \tilde{\th}_i=\pi\}}
\mbox{\hspace{-2ex}}\tilde{\th}_j^2 =\#\{j\in J': \l'_{\b'(j)}=-1\}\pi^2
=\dim(E_{-1}(R_1))\pi^2 = l_+\pi^2,
\ee

\noi and similarly $\sum_{j\in J\setminus J_*} \tilde{\th}_j^2 =  l_-\pi^2.$
From \eqref{dri2}
we therefore have

\bestar
\dsop(RI_\mbs,I)^2 &=&\frac{1}{2}\left\{
\sum_{j\in J'\setminus J'_*} \tilde{\th}_j^2+
\sum_{j\in J\setminus J_*} \tilde{\th}_j^2
+\sum_{j\in J'_*}\tilde{\th}_j^2+\sum_{j\in J_*}\tilde{\th}_j^2\right\}\\
&=&
\frac{1}{2}\left\{
l_+\pi^2+
l_-\pi^2+
2\sum_{j\in J_*}\tilde{\th}_j^2\right\},
\eestar

\noi establishing \eqref{dri3}.

If $\lev(\boldsig)=\lev(R)$, then equation \eqref{lplm2} implies that $l_+=l_-$,
so \eqref{dri3} implies the first equality in \eqref{dri5}.  For the second equality,
observe that $j\in J\setminus J_*$ if and only if $\tilde{\th}_j$ is 0 or $\pi$.
The number of $j$'s in $J$ for which $\tilde{\th}_j=\pi$ is exactly $l_-$, while the $j$'s in $J$ for which $\tilde{\th}_j=0$ have no effect on $\sum_{j\in J}\tilde{\th}_j^2$.  Hence
the second equality in \eqref{dri5} holds.
\qedns

\begin{cor}\label{equivtograss}
Hypotheses and notation as in Lemma
\ref{R-eigenlemma}, except that we additionally write $m_R:=\lev(R)$
and $m=\min\{m_\mbs,m_R\}$.
Let $\{\th_i\in [0,\pi]\}_{i=1}^{\ceil{p/2}}$ be angles for which
$\sfr(\th_1,\dots, \th_{\ceil{p/2}})$ is a normal form of $RI_\mbs$, let
$\{\tilde{\th}_i\}_{i=1}^p$ be as defined in \eqref{deftth}, let
the elements of $J$ be $i_1<i_2<\dots <i_{m_\mbs}$, and let
$\phi_{J,j}=\phi_{J,j}(E_{-1}(R)), 1\leq j\leq m$. Then:

(i)  Up to ordering,

\be\label{relatephilam2}
\phi_{J,j} = \frac{\tilde{\th}_{i_j}}{2}, \ \ 1\leq j\leq m.
\ee

(ii)  If $m_\mbs =m_R$ then
\be\label{dsop2dgr}
\dsop(RI_\mbs,I)=2d_{Gr}(E_{-1}(R),\bfr^J).
\ee
\end{cor}

\pf (i). Let $\Tilde{W}$ be the $p\times m_R$ matrix formed by the columns of the basis \eqref{basisem} of $E_{-1}(R)$, with the elements of the first set in \eqref{basisem} comprising the first $m_\mbs-l_-$ columns , and the elements of
the second set comprising the last $l_+$ columns. (Here $l_\pm$ are defined as in Lemma
\ref{R-eigenlemma}(vii).)
Without loss of generality we order the
$R_4$-eigenvectors $w_i$ such that the first $m_\mbs-l_-$ are the ones for which $\l_i\neq 1$.

Since the
columns of $\Tilde{W}$ form an orthonormal basis of $E_{-1}(R)$, the
numbers $\{\cos\phi_{J,i}\}_{i=1}^m$ are the singular values of the
$m_R\times m_\mbs$ matrix $\Tilde{W}^T {\sf E}_J$.  (This is true whether $m_R\leq m_\mbs$ or $m_R>m_\mbs$.) But, relative to the block-decomposition of matrices used in
Lemma \ref{R-eigenlemma}, the upper $(p-m_\mbs)\times m_\mbs$ block of ${\sf E}_J$ is $0$, and the lower $
m_\mbs\times m_\mbs$ block is $I_{m_\mbs\times m_\mbs}$.  Hence, writing $\Tilde{W}_*$ for the $m_\mbs\times (m_\mbs-l_-)$
matrix formed by the last $m_\mbs$ rows of the first $m_\mbs-l_-$ columns of $\Tilde{W}$, and noting that $m_\mbs-l_-=m_R-l_+$ (by \eqref{lplm2}), we have $\Tilde{W}^T{\sf E}_J = \left[\begin{array}{c} \Tilde{W}_*^T \\ 0_{l_+ \times m_\mbs}\end{array}\right]$, where the $i^{\rm th}$ row of the $(m_\mbs-l_-)\times (p-m_\mbs)$ matrix $\Tilde{W}_*^T$ is a multiple
of $w_i^T$.
Hence for $i,j \leq m_\mbs-l_-=m_R-l_+$,

\be
\left((\Tilde{W}^T{\sf E}_J)(\Tilde{W}^T{\sf E}_J)^T\right)_{ij}=
\sqrt{\frac{1-\l_i}{2}}\sqrt{\frac{1-\l_j}{2}}
\ \lb w_i, w_j\rb
=
\frac{1-\l_j}{2}\d_{ij}\ ,
\ee

\noi and all other entries of the $m_R\times m_R$ matrix  $\Tilde{W}^T{\sf E}_J (\Tilde{W}^T{\sf E}_J)^T$ are 0. But for \linebreak
$m_\mbs-l_- <i\leq m_\mbs$,
we have $\l_i=1$, so $\left((\Tilde{W}^T{\sf E}_J)(\Tilde{W}^T{\sf E}_J)^T\right)_{ij}=
\frac{1-\l_j}{2}\d_{ij}$ for all $i,j \leq m=\min\{m_R,m_\mbs\}$.  Thus the upper left-hand
$m\times m$ block of $(\Tilde{W}^T{\sf E}_J)(\Tilde{W}^T{\sf E}_J)^T$ (the entire $m_R\times m_R$ matrix if $m_R\leq m_\mbs$) is $\diag(\frac{1-\l_1}{2},
\dots,\linebreak \frac{1-\l_m}{2})$, so the numbers $\sqrt{\frac{1-\l_j}{2}}, 1\leq j\leq m,$ are
the singular values of $\Tilde{W}^T{\sf E}_J$.  Thus, up to ordering, the principal
angles $\{\phi_{J,i}\}$ are given by

\be\label{relatephilam}
\cos\phi_{J,j} = \sqrt{\frac{1-\l_j}{2}} \ , 1\leq j\leq m.
\ee

The bijection $\b:J\to \{1,\dots,m_\mbs\}$ used in the proof of Corollary \ref{R-eigencor2} is simply
the inverse of the map $j\mapsto i_j$. Thus from \eqref{relatethlam}, we have

\be\label{relatethlam3}
-\l_j = \cos\tilde{\th}_{i_j}, \ \ 1\leq j\leq m_\mbs\ .
\ee

\noi Combining \eqref{relatephilam} with \eqref{relatethlam3},

\be\label{cosequal}
\cos\phi_{J,j} = \sqrt{\frac{1+\cos\tilde{\th}_{i_j}}{2}}
=\cos\frac{\tilde{\th}_{i_j}}{2}.
\ee

\noi But $\tilde{\th}_{i_j}\in [0,\pi]$, so both $\phi_{J,j}$ and
$\frac{\tilde{\th}_{i_j}}{2}$ lie in $[0,\frac{\pi}{2}]$.  Hence
\eqref{cosequal} implies that $\phi_{J,j}=\tilde{\th}_{i_j}/2$, $1\leq j\leq m$.

\ss
(ii) Assume $m_\mbs=m_R$; then both equal $m$. Corollary \ref{R-eigencor2} then implies that

\be\label{dri4}
\dsop(RI_\mbs,I)^2 =\sum_{i=1}^m  \tilde{\th}_{i_j}^2 .
\ee

\noi But from part (i) we have $\tilde{\th}_{i_j}=2\phi_{J,j}$ for $1\leq j\leq m$, so, using \eqref{dri5},  $\dsop(RI_\mbs,I)^2 =
4\dgr(E_{-1}(R),\bfr^J)^2$, implying \eqref{dsop2dgr}.
\qedns

 We are now ready to establish the general half-angle relation:

\ssn {\bf Proof of Proposition \ref{cor:halfangle}}. Let $U\in O(p)$ and let $T_U:\bfr^p\to \bfr^p$ be the corresponding orthogonal transformation. For any even $m'>0$ and any $R\in {\rm Inv}_{m'}(p)$, we have
\be\label{intertwine}
E_{-1}(URU^{-1}) = T_U(E_{-1}(R)).
\ee

Now let $T:\bfr^p\to\bfr^p$ be an orthogonal transformation carrying
$E_{-1}(R_2)$  to a coordinate plane $\bfr^J$, and let $U\in O(p)$ be the matrix for which $T=T_U$.
Then $UR_2U^{-1}=I_\mbs$, where $\boldsig=\boldsig^J$.
For $i=1,2$ let $R'_i = UR_i U^{-1}$. Since $T$ is an orthogonal transformation, the (multi-)set of principal angles between $E_{-1}(R_1')=T(E_{-1}(R_1))$ and $E_{-1}(R_2')=T(E_{-1}(R_2))$
is identical to the \linebreak (multi-)set of principal angles between $E_{-1}(R_1)$ and $E_{-1}(R_2)$. But
  $R_1'I_\mbs = R_1'R_2'= UR_1R_2U^{-1}$, so $\sfr(\th_1,\dots, \th_{\ceil{p/2}})$ is a normal form of $R_1'I_\mbs$ as well as of $R_1R_2$.  The result now follows from Corollary \ref{equivtograss}(i)
and equation \eqref{relatethlam2} (the latter being needed only for the final statement of the result).
\qedns

\subsection{ The proofs of Propositions
\ref{isometry} and \ref{equivstat}}\label{sec:app.3b}

\noi {\bf Proof of Proposition \ref{isometry}.}
Since $\dsop(RI_\mbs,I)=\dsop(R, I_\mbs^{-1})=
\dsop(R,I_\mbs)$,
conclusion (ii) of  Corollary \ref{equivtograss} can be written equivalently as:
\be\label{almostisom}
\dsop(\Phi(W), \Phi(\bfr^J))=2\dgr(W,\bfr^J).
\ee

Fix any $J\in \jmp$.  Letting ``$\dotprod$'' denote the
natural left-action of $\Sop$ on $\Gr_m(\bfr^p)$,
observe that, in the notation of the proof of Proposition \ref{cor:halfangle}),
 for all $U\in\Sop$ and $W\in \Gr_m(\bfr^p)$
we have $\Phi(U\dotprod W) = U\Phi(W)U^{-1}$ (simply another way of writing \eqref{intertwine}.)  Clearly $\dgr$ is invariant
under this action, and $\dsop$ is both left- and right-invariant, so \eqref{almostisom}
implies that
\ben\dsop(\Phi(U\dotprod W), \Phi(U\dotprod \bfr^J))
=2\dgr(U\dotprod W, U\dotprod \bfr^J).
\een

Now let $W,V\in \Gr_m(\bfr^p)$. Since the action of $\Sop$ on $\Gr_m(\bfr^p)$ is transitive, there exists $U\in \Sop$ such that $U\dotprod \bfr^J=V.$ Using any such $U$, we then have
\bestar
\dsop(\Phi(W), \Phi(V))&=&
\dsop(\Phi(U\dotprod U^{-1}\dotprod W), \Phi(U\dotprod \bfr^J))\\ &=& 2\dgr(U\dotprod U^{-1}\dotprod W, U\dotprod \bfr^J)\\
&=& 2\dgr(W,V).
\eestar
\qedns

\begin{remark}\label{rem:isom} \rm Of course, Proposition \ref{isometry} can
be deduced from computations with the principal fibration $$\pi: \Sop\to
\Sop/S(O(m)\times O(p-m)) \iso \Gr_m(\bfr^p);$$ the standard Riemannian
metric on $\Gr_m(\bfr^p)$ (for which $\dgr$ is the geodesic-distance function) is
defined so as to make $\pi$ a Riemannian submersion up to a normalization constant.
Our proof of Proposition \ref{isometry} is independent of this Riemannian proof in
the sense that it establishes equality between the left-hand side of
\eqref{almostisom} and the right-hand side {\em as defined by equation
\eqref{wongdist}}. Without the {\em a priori} knowledge that $\dgr$ is a
geodesic-distance function, it is not obvious that  $\dgr$ satisfies the triangle
inequality, hence whether $\dgr$ is a metric.  Thus Proposition \ref{isometry}
actually provides an independent proof that $\dgr$ is a metric on $\Gr_m(\bfr^p)$.
The only use of Riemannian geometry in this proof is through the knowledge that
$\dsop$ is, in fact,  a metric (because it is a geodesic-distance function).
\end{remark}

\noi {\bf Proof of Proposition \ref{equivstat}.}
Let ``Statement 1'' and ``Statement 2'' be the statements listed as 1 and 2 in the Proposition.
As noted in the proof of Proposition \ref{isometry},
 $\dsop(RI_\mbs,I)=\dsop(R,I\mbs)$, so the inequality
$\dsop(RI_\mbs, I) <  \dsop(R,I)$ can be rewritten as

\ben
\dsop(R,I_\mbs)^2 <  \frac{m\pi^2}{2} \ .
\een

Assume first that Statement 1 is true.
Let  $W\in \Gr_m(\bfr^p)$. Then $\Phi_{m,p}(W)$
is an involution of level $m$, so there exists  $\boldsig\in\I_p^+$ of level $m$
such that \linebreak $\dsop(\Phi_{m,p}(W),I_\mbs)^2 <  \frac{m\pi^2}{2}$.  Select such a $\boldsig$ and let $J =J^\mbs$.  Then $I_\mbs=\Phi_{m,p}(\bfr^J)$, so

\bestar
\dgr(W,\bfr^J)^2 = \frac{1}{4}\dsop(\Phi_{m,p}(W),\Phi_{m,p}(\bfr^J))^2
&=&\frac{1}{4}\dsop(\Phi_{m,p}(W),I_\mbs)^2\\
& < & \frac{m\pi^2}{8}\ .
\eestar

\noi Hence Statement 2 is true.

Conversely, assume that Statement 2 is true. Let  $R\in {\rm Inv}_m(p)$.  Then there exists $J\in\jmp$ such that
$\dgr(\Phi_{m,p}^{-1}(R), \bfr^J)^2 < \frac{m\pi^2}{8}$. Select such a $J$ and let $\boldsig =\boldsig^J$.  Then $I_\mbs=\Phi_{m,p}(\bfr^J)$, so

\ben
\dsop(R,I_\mbs)^2 = \dsop(R,\Phi_{m,p}(\bfr^J))^2
=4\dgr(\Phi_{m,p}^{-1}(R),\bfr^J)^2\\
<\frac{m\pi^2}{2}\ .
\een

\noi Hence Statement 1 is true.
\qedns

\setcounter{equation}{0}
\section{
Proofs of sign-change reducibility results, part I:
Propositions \ref{involprop-weak} and \ref{conjtruem2}}\label{sec:app.4}

We are now ready to attack the question of sign-change reducibility: given $R\in{\rm Inv}(p)$, can we find
$\boldsig\in \I_p^+$ such that $\dsop(RI_\mbs,I)<\dsop(R,I)$?
 Equations \eqref{invodist2} and \eqref{dri3} tell us that
this inequality is satisfied if and only if

\be\label{goal}
(l_+ + l_-)\pi^2
+2 \sum_{j\in J_*} \tilde{\th}_j^2 < \lev(R)\pi^2,
\ee

\noi where $l_\pm=l_\pm(R,\boldsig)$ are as in Lemma \ref{R-eigenlemma}(vii).  Since $\pi$ is the
largest possible value for a normal-form angle
in \eqref{normform}, it is reasonable to try to look for a $\boldsig$ such that
$l_+$ and $l_-$ are as small as possible.   However, to achieve \eqref{goal},
we have to make sure that we do not make $\sum_{j\in J_*}\tilde{\th}_j^2$
too large while we are making $l_\pm$ small.
We next prove a lemma that, via its subsequent corollary, will help us show that for $\lev(R)=m\geq \frac{p}{2}$,
we can choose $J\in \jmp$ to make $\dsop(RI_{\mbs^J},I)$ as small as is needed to prove Proposition \ref{involprop-weak}.

\begin{lemma}\label{combinat} For $1\leq m\leq p$,

\be\label{eq:combinat}
\sum_{J\in \jmp} {\sf E}_J {\sf E}_J^T = {p-1 \choose m-1} I_{p\times
    p}.
\ee
\end{lemma}

\pf
For $J=(i_1,\dots,i_m)\in\jmp$, we have

\be
{\sf E}_J{\sf E}_J^T = \sum_{i\in J} \bfe_i\bfe_i^T\ .
\ee

\noi Hence when $m=1$ and when $m=p$, the left-hand side of
\eqref{eq:combinat} reduces to $I_{p\times p}$, which is also true of
  the right-hand side.

We proceed by induction on $p$. For each $p\geq 1$, consider the
statement

\be\label{induct}
S(p):\ \mbox{Equation \eqref{eq:combinat} is true for all $m$ satisfying
$1\leq m\leq p$.}
\ee

\noi We have already established that \eqref{eq:combinat} holds for
$m=1=p$, hence that statement $S(1)$ is true.
Now suppose that $S(p)$ is true for some given $p$.  To consider $S(p+1)$, let
$\{\bfe_i\}_{i=1}^{p}, \{\bfe_i'\}_{i=1}^{p+1}$ denote the standard
bases of $\bfr^p, \bfr^{p+1}$ respectively. For $K=\{i_1, \dots,
i_m\}\in{\cal J}_{p+1,m}$ with $i_1<i_2<\dots <i_m$ we write $E'_K$
  for the $(p+1)\times m$ matrix whose $j^{\rm th}$ column is
  $\bfe'_{i_j}$, $1\leq j\leq m$. Note that

\bestar
\bfe_i' &=&
\left[\begin{array}{c}\bfe_i\\ 0\end{array}\right]\ \ \ \mbox{for}\ 1\leq
i\leq p, \\
{\sf E}_J'({\sf E}_J')^T &=& \left[
\begin{array}{cc} {\sf E}_J({\sf E}_J)^T & 0_{p\times 1} \\
0_{1\times p} & 0\end{array}\right]\ \ \ \mbox{for}\ J\in {\cal
  J}_{m,p}, \\
\mbox{and \hspace{.5in}}
\bfe_{p+1}'(\bfe_{p+1}')^T &=& \left[
\begin{array}{cc} 0_{p\times p} & 0_{p\times 1} \\
0_{1\times p} & 1\end{array}\right].
\eestar

\noi Hence for $1\leq m\leq p$,

\bestar
\lefteqn{
\sum_{K\in{\cal J}_{m,p+1}} {\sf E}_K'({\sf E}_K')^T}\\
&=& \sum_{\{K\in{\cal J}_{m,p+1}: p+1\in K\}} {\sf E}_K'({\sf E}_K')^T
+\sum_{\{K\in{\cal J}_{m,p+1}: p+1\notin K\}} {\sf E}_K'({\sf E}_K')^T\\
&=&
\sum_{\{K\in {\cal J}_{m,p+1} :
K=J\cup\{p+1\}\mbox{\scriptsize \ for\ some\ }J\in {\cal J}_{m-1,p}\}}
\mbox{\hspace{-.75in}}{\sf E}_K'({\sf E}_K')^T \mbox{\hspace{.5in}}
+\sum_{K\in\jmp} {\sf E}_K'({\sf E}_K')^T\\
&=&
\sum_{J\in{\cal J}_{m-1,p}}
\left( {\sf E}_J'({\sf E}_J')^T + \bfe_{p+1}'(\bfe_{p+1}')^T\right)
       +\sum_{J\in\jmp} \left[
\begin{array}{cc} {\sf E}_J({\sf E}_J)^T & 0_{p\times 1} \\
0_{1\times p} & 0\end{array}\right]\\
&=&
\sum_{J\in{\cal J}_{m-1,p}}
\left( \left[
\begin{array}{cc}{\sf E}_J{\sf E}_J^T &
  0_{p\times 1}\\
0_{1\times p} & 0 \end{array}\right]\right)
+|{\cal J}_{m-1,p}|\,\bfe_{p+1}'(\bfe_{p+1}')^T\\
&&\mbox{\hspace{.5in}}+\sum_{J\in\jmp} \left[
\begin{array}{cc} {\sf E}_J({\sf E}_J)^T & 0_{p\times 1} \\
0_{1\times p} & 0\end{array}\right]\\
&=&
\left[
\begin{array}{cc}\sum_{J\in{\cal J}_{m-1,p}} {\sf E}_J{\sf E}_J^T
+\sum_{J\in\jmp} {\sf E}_J{\sf E}_J^T &
  0_{p\times 1}\\
0_{1\times p} & 0 \end{array}\right]
+|{\cal J}_{m-1,p}|\,\bfe_{p+1}'(\bfe_{p+1}')^T\\
&=&
\left[
\begin{array}{cc}\left\{{p-1\choose m-2} + {p-1\choose m-1}\right\}I_{p\times p} &
  0_{p\times 1}\\
0_{1\times p} & 0 \end{array}\right]
+{p\choose m-1}\, \bfe_{p+1}'(\bfe_{p+1}')^T
\\
&=&{p\choose m-1}I_{(p+1)\times (p+1)}\ .
\eestar

Hence \eqref{combinat} holds with $p$ replaced by $p+1$, as long as
$1\leq m\leq p$.  But we have already established that
\eqref{combinat} holds whenever $m=p$; hence if $p$ is replaced by
$p+1$, the equality holds for $m=p+1$.  Thus \eqref{combinat} holds
for all $m$ with $1\leq m\leq p+1$; i.e. statement $S(p+1)$ is true.
By induction, $S(p)$ is true for all $p$, which is exactly what the
Lemma asserts.\qedns

\begin{cor}\label{princanglebound}
Let $m\in \{1,2,\dots,p\}$ and let $W\in \Gr_m(\bfr^p)$. There exists
$J\in \jmp$ such that

\be\label{eq:princanglebound}
\sum_{i=1}^m \sin^2 \phi_{J,i} \leq m(1-\frac{m}{p}).
\ee

\noi Furthermore, the inequality in \eqref{eq:princanglebound} is
strict for {\em some} $J\in\jmp$ unless equality holds in
\eqref{eq:princanglebound} for {\em all} $J\in\jmp$.

\end{cor}

\pf Let $\Tilde{W}$ be any $p\times m$ matrix whose columns are an
orthonormal basis of $W$.  Using Lemma \ref{combinat},

\bestar
\sum_{J\in \jmp} \tr(\Tilde{W}^T {\sf E}_J {\sf E}_J^T\Tilde{W})
&=& \tr\left(\Tilde{W}^T\left(\sum_{J\in \jmp} {\sf E}_J {\sf E}_J^T\right)\Tilde{W}\right)\\
&=& \tr\left(\Tilde{W}^T{p-1 \choose m-1}I_{p\times p}\Tilde{W}\right)\\
&=& m{p-1 \choose m-1}
\eestar

\noi since $\Tilde{W}^T\Tilde{W}=I_{m\times m}$.

 Since $|\jmp| = {p \choose m}$, the average of $\tr(\Tilde{W}^T {\sf E}_J
{\sf E}_J^T\Tilde{W})$ over all $J\in \jmp$ is $m{p-1 \choose m-1}/{p
  \choose m} = m^2/p$.  Hence $\tr(\Tilde{W}^T {\sf E}_J {\sf E}_J^T\Tilde{W})\geq
m^2/p$ for at least one $J\in \jmp$, and the inequality is strict for
some $J$ unless it is an equality for all $J$.  But for any $Z\in
\Gr_m(\bfr^p)$, the principal angles $\phi_1, \dots, \phi_m$ between
$W$ and $Z$ are the numbers in $[0,\frac{\pi}{2}]$ for which
$\cos\phi_1, \dots, \cos\phi_m$ are the singular values of the
$m\times m$ matrix $\Tilde{W}^T\Tilde{Z}$, where $\Tilde{Z}$ is any
$p\times m$ matrix whose columns are an orthonormal basis of $Z$.
Since for any $J\in \jmp$ the columns of ${\sf E}_J$ are an orthnormal basis
of $\bfr^J$, it follows that $\sum_{i=1}^m
\cos^2\phi_{J,i}=\tr(\Tilde{W}^T {\sf E}_J (\Tilde{W}^T {\sf E}_J )^T)
=\tr(\Tilde{W}^T {\sf E}_J {\sf E}_J^T\Tilde{W})$. Thus, for some $J$,
$\sum_{i=1}^m \cos^2 \phi_{J,i} \geq \frac{m^2}{p}$, and the
inequality is strict for some $J$ unless it is an equality for all
$J$.  But for any given $J$,

\be\label{sin2bd}
\sum_{i=1}^m \cos^2 \phi_{J,i} \geq
\frac{m^2}{p}\iff
\sum_{i=1}^m\sin^2\phi_{J,i} = m-\sum_{i=1}^m
\cos^2\phi_{J,i} \leq m-\frac{m^2}{p} = m(1-\frac{m}{p}) \ ,
\ee

\noi and the first inequality in \eqref{sin2bd} is strict if and only if the second is strict. Thus \eqref{eq:princanglebound} holds for some $J$, and the
inequality in \eqref{eq:princanglebound} is strict for some $J$ unless it is an equality for all
$J$.
\qedns

\noi{\bf Proof of Proposition \ref{involprop-weak}}.

If $m=p$  then $p$ is even, $R=-I$, and for
$\boldsig=(-1,-1,\dots, -1)$ we have $I_\mbs=-I$ and
$\dsop(RI_\mbs,I)=0<\dsop(R,I)$. Henceforth we assume $m<p$.

Let $W=E_{-1}(R)$ and let $m=\dim(W)$.  Note that

\be\label{dinvi}
\dsop(R,I)^2 = \frac{m}{2}\pi^2.
\ee

Let $J\in\jmp$ be such that $\sum_{i=1}^m \sin^2 \phi_{J,i} =
\min_{K\in\jmp}\{\sum_{i=1}^m \sin^2 \phi_{K,i}\}$.
By Corollary \ref{princanglebound}, inequality
\eqref{eq:princanglebound} holds, and the inequality is strict unless

\be\label{sin2bdeq}
\sum_{i=1}^m \sin^2 \phi_{K,i}=m(1-\frac{m}{p})
\ee

\noi  for all $K\in\jmp$.
Let $\boldsig=\boldsig^J$.
By Corollary \ref{equivtograss},

\be\label{dist2eq}
\dsop(RI_\mbs)^2 =4 d_{Gr}(W,\bfr^J)^2=4\sum_{i=1}^m (\phi_{J,i})^2
\ee

\noi where $\phi_{J,i}=\phi_{J,i}(W)$.

The function $f: x\mapsto \frac{\sin x}{x}$ is
strictly decreasing on the interval $(0,\frac{\pi}{2}]$.  Hence for
all $x\in(0,\frac{\pi}{2}]$ we have $\frac{\sin x}{x}\geq
f(\frac{\pi}{2})=\frac{2}{\pi}$, with equality only if
$x=\frac{\pi}{2}$; thus for $x\in [0,\frac{\pi}{2}]$ we have $x\leq \frac{\pi}{2}\sin x$, with equality
only if $x=0$ or $x=\frac{\pi}{2}$.
Hence

\bearray
\label{ineq1}
\dsop(RI_\mbs)^2 = 4\sum_{i=1}^m (\phi_{J,i})^2 &\leq&
\pi^2\sum_{i=1}^m \sin^2\phi_{J,i}\\
\label{ineq2}
&\leq & m(1-\frac{m}{p})\pi^2\\
\label{ineq3}
&\leq & \frac{m}{2}\pi^2 \ \ \ \mbox{(since $\frac{m}{p}\geq \frac{1}{2})$}\\
\nonumber
&=& \dsop(R,I)^2.
\eearray

Hence $\dsop(RI_\mbs)\leq\dsop(R,I)$, and this inequality is strict if
any of the inequalities \eqref{ineq1}, \eqref{ineq2}, \eqref{ineq3}
is strict.  Inequality \eqref{ineq1} is strict if $0<
\phi_{J,i}<\frac{\pi}{2}$ for some $i$, and, by our choice of $J$,
\eqref{ineq2} is strict unless equality holds in \eqref{sin2bdeq}
for all $K\in\jmp$.

We claim that at least one of the inequalities \eqref{ineq1}, \eqref{ineq2}
is strict.  Assume this is not so.  Then, since equality holds in \eqref{ineq1}
with $J$ replaced by any $K\in\jmp,$ it follows that for {\em all} $K\in\jmp$ and $i\in\{1,\dots,m\}$ the angle $\phi_{K,i}$ is either 0 or $\pi/2$, and
that $\sum_{i=1}^m \sin^2 \phi_{K,i} = m(1-\frac{m}{p})$ for all $K$.  But for any $V\in \Gr_m(\bfr^p)$, there always exists $K\in\jmp$ for which none of the
principal angles $\phi_{K,i}(V,\bfr^{K})$ is $\pi/2$. Choosing such $K$
for our $m$-plane $W,$ all of the principal angles $\phi_{K,i}$ must therefore be 0 (since they are all either 0 or $\pi/2$).  But then
$\sum_{i=1}^m \sin^2 \phi_{K,i} =0 < m(1-\frac{m}{p})$, a contradiction.

Thus at least one of the inequalities \eqref{ineq1}, \eqref{ineq2} is strict,
so $\dsop(RI_\mbs)<\dsop(R,I)$.
\qedns

We will establish
Proposition \ref{conjtruem2} (a weak version of Conjecture \ref{onlymcanwork})
as a consequence of a different weakened version of  Conjecture \ref{onlymcanwork}:

\begin{prop}\label{weakconjtrue}
Let $m\geq 2$ be even, let $R\in \Sop$ be an involution of level $m$, and let $\boldsig\in \I_p^+$.  If $\dsop(RI_\mbs,I)<\dsop(R,I)$, then $\lev(\boldsig)<2m$.
(Hence if $R$ is sign-change reducible, then it is reducible by a sign-change of level
less than $2m$.)
\end{prop}

This proposition, which we will prove this shortly,
reduces Proposition \ref{conjtruem2} into a triviality:

\ssn {\bf Proof of Proposition \ref{conjtruem2}, assuming Proposition \ref{weakconjtrue}:}
The only positive even integer less than $2\times 2$ is 2.\qedns

\noi{\bf Proof of Proposition \ref{weakconjtrue}}.
Let $m_\mbs=\lev(\boldsig)$. Define $l_\pm$
as in Lemma \ref{R-eigenlemma}.
From \eqref{dri3},

\ben
\frac{1}{2}(l_+ + l_ -) \pi^2 \leq \dsop(RI\mbs,I)^2 < \dsop(R,I)^2
=\frac{1}{2}m_R\pi^2,
\een

\noi so

\be\label{lppluslm}
l_+ + l_- < m_R\ .
\ee

\noi But by \eqref{lplm2} we have $l_- = l_+ + m_\mbs - m_R$, so substituting
into  \eqref{lppluslm}, we have $2l_+ +m_\mbs -m_R <m_R$; equivalently,

\ben
2l_+ <2m_R-m_\mbs .
\een

\noi Since $l_+\geq 0$, we must have $m_\mbs<2m_R$. \qedns

\setcounter{equation}{0}
\section{Proofs of sign-change reducibility results, part II: Proposition  \ref{pleq4prop}
}\label{sec:app.5}

As noted in Section \ref{sec:app.2}, Proposition \ref{involprop-weak} proves part (a) of Proposition \ref{pleq4prop}.  Thus it remains only to prove part (b) of this Proposition.

The combination of Proposition \ref{conjtruem2} and Proposition \ref{equivstat} is what will guide our proof of
part (b). To
establish the result, it suffices to prove that for $p\geq 11$, the answer to Question \ref{involprop} is no---i.e. that there exist involutions in $\Sop$ that are not sign-change reducible. Hence it suffices to prove that there exist such involutions of level 2.  By Proposition \ref{conjtruem2}, it therefore suffices to establish (for $p\geq 11$) the existence of involutions that are not sign-change reducible by a sign-change of level 2; thus it suffices to  show that Statement 2 of Proposition \ref{equivstat} is false when $p\geq 11$ and $m=2$.  For this, we need only produce planes in $\bfr^p$ for which we can show that
\eqref{grassineq} is false for all $J\in {\cal J}_{2,p}$. Towards this end, we examine two (families) of examples  in which $m=2$ and $p\geq 4$.

\begin{Example}\label{counterex} {\rm Let $p=2k$ or $2k+1$, where $k\geq 2$.
Define vectors $\hat{v},\hat{w}\in\bfr^p$ by

\ben
\hat{v}  = \frac{1}{\sqrt{k}}\sum_{i=1}^k \bfe_i, \ \ \ \ \
\hat{w} = \frac{1}{\sqrt{k}}\sum_{i=k+1}^{2k}\bfe_i\ .
\een

The set $\{\hat{v},\hat{w}\}$ is orthonormal.  Let $W_p=\spn\{\hat{v},\hat{w}\}$,
a 2-plane in $\bfr^p$. We will compute the principal angles between $W_p$ and
$\bfr^J$ for all $J\in {\cal J}_{2,p}$. Write $J=\{i,j\}$, where $1\leq i<j\leq p$.
Let $\Tilde{W}$ be the $p\times 2$ matrix whose first column is $\hat{v}$ and whose second column is $\hat{w}$.  Since the columns of $\Tilde{W}$ are an orthonormal basis of
$W_p$, the principal angles between $W_p$ and $\bfr^J$ are the arc-cosines
of the singular values of $\Tilde{W}^T {\sf E}_J$.

First suppose that $p$ is even.  We divide the elements $\{i,j\}\in {\cal J}_{2,p}$ into two cases:
Case I= $\{\{i,j\} : i<j\leq k\ \mbox{ or}\  k< i<j\}$;
Case II= $\{\{i,j\} : i\leq k<j\}$.  The principal values of the $2\times 2$ matrix $\Tilde{W}^T {\sf E}_J$
are easily computed to be $0$ and $\frac{4}{p}$ in Case I, and $\frac{2}{p}$ (with multiplicity 2)
in Case II.  Hence the principal angles are

\begin{align*}
\phi_{J,1} =\frac{\pi}{2}, \ \ & \phi_{J,2} =\cos^{-1}\sqrt{4/p}
 \ \ \mbox{in Case I}, \\
\phi_{J,1} =\phi_{J,2}  &=\cos^{-1}\sqrt{2/p}
 \ \ \mbox{in Case II},
\end{align*}

\noi so

\be\label{peven}
\min_{J\in {\cal J}_{2,p}}\{\dgr(W_p,\bfr^J)^2\}
=\min\left\{\left(\frac{\pi}{2}\right)^2 + \left(\cos^{-1}\sqrt{4/p}\right)^2,
2\left(\cos^{-1}\sqrt{2/p}\right)^2\right\}.
\ee

We will return to \eqref{peven} shortly, but first let us do the analogous computation
for $p$ odd.  For $p=2k+1$, we divide the computation into three cases:
Case I= $\{\{i,j\} : i<j\leq k\ \mbox{ or}\  k< i<j\leq 2k\}$;
Case II= $\{\{i,j\} : i\leq k<j\leq 2k\}$; and Case III= $\{\{i,j\} : i\leq 2k, j=2k+1\}$.
The principal values of the matrix
$\Tilde{W}^T {\sf E}_J$ are $0$ and $\frac{2}{p}+\frac{2}{p-1}$
in Case I, $\frac{2}{p}$ and $\frac{2}{p-1}$ in Case II, and
$\frac{2}{p-1}$ in Case III.
Hence the principal angles are

\begin{align*}
\phi_{J,1} =\frac{\pi}{2}, \ \ & \phi_{J,2} =\cos^{-1}\sqrt{4/(p-1)}
 \ \ \mbox{in Case I}, \\
\phi_{J,1} =\phi_{J,2}  &=\cos^{-1}\sqrt{2/(p-1)}
 \ \ \mbox{in Case II},\\
\phi_{J,1} =\frac{\pi}{2}, \ \ & \phi_{J,2} =\cos^{-1}\sqrt{2/(p-1)}
 \ \ \mbox{in Case III}.
\end{align*}

\noi Clearly $\phi_{J,1}^2 +\phi_{J,2}^2$ is larger in Case III than in Case II,
so

\be\label{podd}
\min_{J\in {\cal J}_{2,p}}\{\dgr(W_p,\bfr^J)^2\}
=
\min\left\{\left(\frac{\pi}{2}\right)^2 + \left(\cos^{-1}\sqrt{\frac{4}{p-1}}\right)^2,
2\left(\cos^{-1}\sqrt{\frac{2}{p-1}}\right)^2
\right\}.
\ee
}
\end{Example}

It follows from \eqref{peven} and \eqref{podd} that

\bearray\label{counterex1}
\lim_{p\to\infty}
\min_{J\in {\cal J}_{2,p}}\{\dgr(W_p,\bfr^J)^2\}=
2 \left(\frac{\pi}{2}\right)^2 &=&\frac{\pi^2}{2}
\\
&\not < &\frac{\pi^2}{4} = \frac{m\pi^2}{8}\
\nonumber
\eearray

\noi since $m=2$ in Example \ref{counterex}. Hence for large enough $p$, Statement 2 in Proposition \ref{equivstat} is false,
and therefore so is Statement 1.  This already shows that for all $p$ sufficiently large,
there exist geodesically antipodal pairs $(U,V)$ in $\Sop\times \Sop$ that are not sign-change
reducible.  However, to get the quantitative statement in Proposition \ref{pleq4prop}(b), we have to
continue working.

It can be shown\footnote{The authors did not find this exercise in Calculus 1 entirely trivial,
but are nonetheless leaving it to the reader.}
 that for $0<x\leq 1$,

\be\label{compar1}
(\pi/2)^2 + (\cos^{-1}x)^2 > 2(\cos^{-1}\frac{x}{\sqrt{2}})^2,
\ee

\noi hence that in  \eqref{peven} in
\eqref{podd}, the second of the two expressions being compared is the smaller. Thus

\be\label{counterex_dist}
\min_{J\in {\cal J}_{2,p}}\{\dgr(W_p,\bfr^J)\}
= \sqrt{2}\cos^{-1}(c_p), \ \ \mbox{where}\ c_p=
\left\{\begin{array}{ll} \sqrt{2/p}, & p\ \mbox{even},\\
\sqrt{2/(p-1)}, & p\ \mbox{odd}.
\end{array}\right.
\ee

Since $m=2$ in Example \ref{counterex}, $\sqrt{m\pi^2/8}=\frac{\pi}{2}$, so equation \eqref{counterex_dist} shows that \eqref{grassineq} (with $W=W_p$)
is
false for all $J\in {\cal J}_{2,p}$ if $\sqrt{2}\cos^{-1}(c_p)\geq \frac{\pi}{2}$;
equivalently, if $c_p\leq \cos\frac{\pi}{2\sqrt{2}} \approx 0.4440.$
This translates to $2\floor{\frac{p}{2}}\geq 2\sec^2\frac{\pi}{2\sqrt{2}}  \approx 10.14$. Hence
the answer to Question \ref{involprop} is definitely ``no" for all $p\geq 12$.  To complete the proof of
Proposition \ref{pleq4prop}(b), it remains only to show that this ``12" can be reduced to ``11".
We will accomplish this with the next example.

\begin{Example}\label{counterex-2}
{\rm Let $p=2k+1$, where $k\geq 2.$
Define vectors $v,w,\hat{v}, \hat{w}\in\bfr^p$ by

\begin{align*}
v  &= \sum_{i=1}^p \bfe_i, &
w &= \sum_{i=1}^k \bfe_i - \sum_{i=k+1}^{2k}\bfe_i\\
\hat{v} &= \frac{v}{\|v\|}=\frac{1}{\sqrt{p}} v, &
\hat{w}&= \frac{w}{\|w\|} = \frac{w}{\sqrt{p-1}} \ .
\end{align*}

As in the previous example, $\{\hat{v},\hat{w}\}$ is an orthonormal basis of a plane $W_p'$.
Just as in Example \ref{counterex}, we can compute the principal angles between $W_p'$ and
$\bfr^J$ for all $J\in {\cal J}_{2,p}$.   We define Cases I and II
and III just as in the odd-$p$ case of the previous example. The principal values
of the relevant $2\times 2$ matrices
are $0$ and $\frac{2}{p}+2/(p-1)$
in Case I, $\frac{2}{p}$ and $\frac{2}{p-1}$ in Case II, and

\ben
\l_\pm(p) := \frac{1}{p}+\frac{1}{2(p-1)} \pm
\sqrt{\frac{1}{p^2}+\frac{1}{4(p-1)^2}}
\een

\noi in Case III. Hence

\bearray\nonumber
\min_{J\in {\cal J}_{2,p}}\{\dgr(W_p',\bfr^J)^2\}
&=&
\min\left\{\left(\frac{\pi}{2}\right)^2 + \left(\cos^{-1}\sqrt{\frac{2}{p}+\frac{2}{p-1}}\right)^2, \right.\\
\nonumber &&\phantom{minx}
\left(\cos^{-1}\sqrt{2/p}\right)^2 +
\left(\cos^{-1}\sqrt{2/(p-1)}\right)^2, \\
&&\phantom{minx}\left.
\left(\cos^{-1}(\sqrt{\l_+(p)})\right)^2 + \left(\cos^{-1}(\sqrt{\l_-(p)})\right)^2
\right\}
\nonumber
\\
\label{counterex-2_distsq}
\eearray

}
\end{Example}

Numerically, we find that
for $p=11$, the middle line of \eqref{counterex-2_distsq} is the smallest of the three lines, so

\bearray\nonumber\label{counterex-2_distsq2}
\min_{J\in {\cal J}_{2,11}}\{\dgr(W_{11}',\bfr^J)^2\} &=&
\left(\cos^{-1}\sqrt{2/11}\right)^2 +
\left(\cos^{-1}\sqrt{2/10}\right)^2\\
&\approx &1.0146\, \frac{\pi^2}{4}\ .
\eearray

Since this number is  larger than $\frac{\pi^2}{4}$, the answer to Question \ref{involprop} is no for $p=11$. This completes the proof of Proposition \ref{pleq4prop}. \qedns

\begin{remarks}\rm (1)  We considered Example \ref{counterex-2} only for odd $p$
because for even $p$, the principal angles $\phi_{J,i}(W_p')$ turn out to be the same
as for $\phi_{J,i}(W_p)$ in Example \ref{counterex}.  In Example \ref{counterex-2},
we can also compute numerically that for $p=5, 7$, and $9$, we have $\min_{J\in {\cal J}_{2,p}}\{\dgr(W_p',\bfr^J)^2\}<\frac{\pi^2}{4}$. However, we cannot conclude
that the answer to Question \ref{involprop} is ``yes" for $p\leq 10$,  since we have
not proven that this example represents the worst case, i.e. that $\min_{J\in {\cal J}_{2,p}}\{\dgr(W_p',\bfr^J)\}\geq \min_{J\in {\cal J}_{2,p}}\{\dgr(W,\bfr^J)\}$
for all $W\in \Gr_m(\bfr^p)$.  Thus Question \ref{involprop} remains open for
$5\leq p\leq 10$.  However, based on computations, it seems likely to the authors
that the largest $p$ for which the answer to Question \ref{involprop} is yes is
closer to 10 than to 4.

(2)
The number $\frac{\pi^2}{2}$ in\eqref{counterex1} is exactly the squared diameter of
$\Gr_2(\bfr^p)$ for all $p\geq 4$.  Thus, \eqref{counterex1} shows that as $p\to\infty$,
the distance between $W_p$ and the {\em closest} coordinate plane(s) $\bfr^J$ is approaching the {\em largest} possible distance between two points in $\Gr_2(\bfr^p)$.
\end{remarks}

\setcounter{equation}{0}
\setcounter{theorem}{0}
\renewcommand{\thesection}{Appendix A}

\section{Partitions and Fibers}\label{sec:fibers}
\renewcommand{\thesection}{A}
\renewcommand{\theequation}{A.\arabic{equation}}

\subsection{Partitions and eigenstructure}
\label{sec:partn}

The strata of each of the stratified spaces in this paper are labeled naturally
either by $\partpset$ or by $\partp.$

The natural left-action of the symmetric group $S_p$ on $\{1,2, \dots, p\}$
induces left-actions of $S_p$ on $\partpset$ and $\bfr^p$.  There is a
canonical bijection between the quotient $\partpset/S_p$ and the set
$\partp$, so we implicitly regard these as the same set. For $\J\in\partpset$,
we write $[\J]$ for the image of $\J$ in $\partp$ under the quotient map.

The sets $\partpset$ and $\partp$ are partially ordered by the refinement
relation. For $\J,\K \in\partpset$, we write $\J\leq \K$ if $\K$ refines $\J$.
Similarly, for $[\J], [\K]\in\partp$ we write $[\J]\leq [\K]$ if $[\K]$
refines $[\J]$. In each of these partially ordered sets there is a well-defined
``highest" (most refined) and ``lowest" (least refined) element; we denote
these with the subscripts ``top" and ``bot" respectively.

\begin{notation}\label{defjd}
{\rm
1. For $D=\diag(d_1,\dots, d_p)\in {\rm Diag}(p)$, let $\J_D$ denote
the partition of $\{1, 2, \dots, p\}$ determined by the
equivalence relation $i\sim_D j \iff d_i=d_j$.

2. For $\emptyset \neq J \subset \{1, 2, \dots, p\}$, let $\bfr^J\subset
\bfr^p$ denote the subspace \linebreak $\{(x_1,\dots, x_p)\in
\bfr^p\mid x_j = 0 \ \forall j\notin J\}$. For a partition $\J=\{J_1,
\dots, J_r\}$ of $\{1, 2, \dots, p\}$ (where the $J_i$ are the blocks of $\J$),
let $\{W_1, \dots, W_r\} =\{W_1^\J, \dots, W_r^\J\}
\linebreak
= \{\bfr^{J_1}, \dots,
\bfr^{J_r}\}$ denote the corresponding subspaces of $\bfr^p$; note
that we have an orthogonal decomposition $\bfr^p = \bfr^{J_1}\plus
\dots \plus \bfr^{J_r}$. Define the subgroup $G_\J\subset\Sop$
by
\be\label{blockstab}
G_\J=\{R\in \Sop\mid RW_i = W_i, 1\leq i\leq r\}.
\ee
\noi
We write $G_\J^0$ for the identity component of $G_\J$.}
\end{notation}

As the reader may check, the above definition of $G_\J$ agrees with the
definition in Section \ref{sec:notn_prelims}:  for all $D\in \D(p)$ we have $G_D=G_{\J_D}$.

For any subgroup $H\subset O(p)$, we write $S(H)$ for $H\intersect \So(p)$.
Note that
\bearray\label{gjiso}
G_\J &\iso& S(O(W_1)\times O(W_2)\times \dots O(W_r))\\
&\iso& S(O(|J_1|)\times O(|J_2|) \dots \times O(|J_r|)),
\label{gjiso2}
\eearray
\noi  where $O(W_i)$ denotes the orthogonal group of the subspace
$W_i$, which we identify with a subgroup of $O(p)$.
Hence, writing $k_i=|J_i|$, we have
\be\label{gj0iso}
G_\J^0\iso SO(k_1)\times SO(k_2) \dots \times SO(k_r).
\ee

\subsection{Signed permutations and signed-permutation matrices}
\label{gpthy}

Let $\I_p =(\ztwo)^p$. The role of $\ztwo$ will be as the {\em
  group of signs}, so we write  its elements as $\pm 1$.
We write the identity element of $\I_p$ as $\bone$.
For $\e\in\ztwo$ and $\boldsig=(\sigma_1,\sigma_2\dots, \sigma_p)\in \I_p$ we
define $\e\boldsig=(\e\sigma_1,\e\sigma_2,\dots, \e\sigma_p)$.

Both $\I_p$ and $S_p$ have natural representations on $\bfr^p$ via sign-changes
and permutations of coordinates, respectively. These representations, which we
denote respectively as $\boldsig\mapsto I_\mbs$ and $\pi\mapsto P_\pi$,
embed $\I_p$ and $S_p$ as subgroups of $O(p)$, together generating
the group of ``signed-permutation matrices''.  Abstractly, this group is a semidirect product $\tsp = \I_p\rtimes S_p$, a split extension of $S_p$ by $\I_p$,
embedded naturally in $O(p)$ via $(\boldsig,\pi)\mapsto I_\mbs P_\pi$.
Defining homomorphisms $\sgn:\I_p\to \ztwo$ and $\Tilde{\sgn}:\tsp\to\ztwo$
 by $\sgn(\sigma_1,\dots,\sigma_p)=\prod_{i=1}^p \sigma_i$ and  $\Tilde{\sgn}(\boldsig,\pi)=\sgn(\boldsig)\sgn(\pi)$ (where $\sgn(\pi)$ is the sign of the permutation $\pi$),
we have $\Tilde{\sgn}(\boldsig,\pi)=\det(I_\mbs P_\pi)$.
Thus the group $\tsp^+$ of even signed-permutations, defined in Section
\ref{sec:notn_prelims}, is simply the kernel of $\Tilde{\sgn}$, and
we have a short exact sequence
\be\label{ses2}
1\ \to\ \I_p^+\ \stackrel{\incl}{\longrightarrow}\ \tsp^+
\ \stackrel{\proj_2}{\longrightarrow}\ S_p\ \to\ 1.
\ee
\noi Since $\I_p^+\iso (\ztwo)^{p-1}$ (non-canonically), $\tsp^+$ is an extension of $S_p$ by $(\ztwo)^{p-1}$, and $|\tsp^+| = 2^{p-1}p!$.

The group $\tsp$ is a well-studied group encountered in other settings (rather different
from this paper's) as $W(B_p)$, the Weyl group of the simple Lie algebra $B_p=\so(2p+1,\bfc)$
\cite{samelson}. Thus $\tsp^+$ is an index-two subgroup of $W(B_p)$. The application
to eigenstructure motivates viewing $\tsp^+$ as an extension of $S_p$: an element of
$\sympp$ determines an element of $\D^+(p)$ up to the action of $S_p$, but this action
does not lift canonically to a fiber-preserving action of $S_p$ on $M(p)$ (at least not for $p$
even; see below); we need to extend $S_p$ to a larger group to obtain such an action.
For each $X\in \sympp$, the fiber $\E_X$ can be identified with positively oriented orthonormal
$X$-eigenbases of $\bfr^p$; the action of $\tsp^+$ sends one such $X$-eigenbasis to another.

A familiar index-two subgroup of $\tsp$ different from $\tsp^+$ is the kernel of the map
$(\boldsig,\pi)\mapsto \sgn(\boldsig)$.  For $p\geq 4$, the latter subgroup is the Weyl group
$W(D_p)$ of the simple Lie algebra $D_p=\so(2p,\bfc)$.  However, the analog of \eqref{ses2}
for $W(D_p)$ splits for all $p$, while \eqref{ses2} splits if and only if $p$ is odd.
For $p$ odd, the map $\tsp^+\to W(D_p)$ defined by $(\boldsig,\pi)\mapsto
(\sgn(\boldsig)\boldsig,\pi)$ is an isomorphism, but it is known that for $p$ even,
$\tsp^+$ is not isomorphic to $W(D_p)$ \cite[p. 151]{pahlings}.

\begin{remark}\label{rem:gjcomps}{\rm For a subspace $W\subset
\bfr^p$ and $\e\in\ztwo=\{\pm 1\}$, let $O_\e(W)\subset O(W)$ denote
the set of orthogonal transformations with determinant $\e$. In the
setting of \eqref{gjiso}, the connected components of $G_\J$ are
$O_{\e_1}(W_1)\times O_{\e_2}(W_2)\times \dots \times O_{\e_r}(W_r)$,
subject to the restriction $\prod_i\e_i=1$.  Thus a labeling of the
blocks of an $r$-block partition $\J$ yields a 1-1 correspondence
between $\I_r^+$ and the set of connected components of $G_\J$. In
particular, the number of connected components is $2^{r-1}$.  }
\end{remark}

Identifying $\D(p)$ with $\bfr^p$, the natural left-action of $S_p$ on $\bfr^p$ yields
a left-action of $S_p$ on $\D(p)$. For $D\in\D(p)$, we will write $[D]$ for its image
in the quotient space $\D(p)/S_p$.

Note that the action of $S_p$ on $\Dpp\subset \D(p)$ lifts to an action of $\tsp$
on $\Dpp$,
\be\label{defgdotd}
g\dotprod D : = \pi_g\dotprod D.
\ee
It is easily seen that $P_g D P_g^{-1} = \pi_g\dotprod D$ for all $g\in\tsp, D\in\D(p)$.

\subsection{Structure of the fibers}\label{fibstr}

The starting point for a systematic description of the fibers of $F$
 is the following proposition. The group-action notation is as in \eqref{act1}.

\begin{prop}\label{jsg1thm3.3} Let $X\in\sympp$ and $(U,D)\in
\E_X=F^{-1}(X)$.
Then
\be\label{fibform1}
\E_X= \{ g\dotprod(UR, D) : R\in G_D,g\in\tsp^+\}.
\ee
\end{prop}

\pf This is a simple corollary of \cite[Theorem 3.3]{JSG2015scarot}.  Details are left to the reader.
\qedns

\begin{cor}\label{cor:fibform4}
Let $X\in\sympp$ and $(U,D)\in \E_X$. Then
\be\label{fibform4}
\E_X= \{ g\dotprod(UR, D) : R\in G_D^0, g\in \tsp^+\}.
\ee
\end{cor}

\pf Clearly the right-hand side of \eqref{fibform4} is contained in
the right-hand side of \eqref{fibform1}, so it suffices to prove the
opposite inclusion.

Let $R\in G_D, g\in\tsp^+$.  Enumerate the blocks of $\J:=\J_D$ as
$J_1, \dots, J_r,$ and let $W_i$ be as in Notation \ref{defjd}. As
noted in Remark \ref{rem:gjcomps}, the enumeration of the blocks of
$\J$ yields a 1-1 correspondence between $\I_r^+$ and the connected
components of $G_\J$.  Let $R$ lie in the component of $G_\J$ labeled by
$(\e_1,\dots, \e_r)\in\I_r^+$. The cardinality of $\{j: \e_j=-1\}$ is
some even number $k$.  Let $\boldsig=(\sigma_1,\dots,\sigma_p)\in
\I_p$, where for $1\leq i\leq p$ we set
\ben
\sigma_i=\left\{\begin{array}{rl}
-1 & \mbox{if}\ i\in J_j,\ \e_j=-1,\ \mbox{and $i$ is the smallest
  element of $J_j$}\ ;\\
1 & \mbox{otherwise}.\end{array}\right.
\een
\noi Then $R_1:=I_\mbs R \in G_\J^0$. But also
$\left|\{i: \sigma_i=-1\}\right|=k$, so $\boldsig\in \I_p^+\subset \tsp^+$, and
$P_g I_\mbs=P_{g_1}$ for some $g_1\in \tsp^+$ with $\pi_{g_1}=\pi_g$. Hence
$P_gR=(P_gI_\mbs)(I_\mbs R)=P_{g_1} R_1$, so
$(U(P_gR)^{-1}, \pi_g\dotprod D) = (U(P_{g_1}R_1)^{-1}, \pi_{g_1}\dotprod D),$
which lies in the right-hand side of \eqref{fibform4}. The
desired inclusion follows. \qedns

To complete our characterization of the fibers of $F$, we introduce one more bit of
notation:

\begin{notation}
\rm
For $\J=\{J_1,\dots, J_r\}\in\partpset$, define
\be
\I_\J^+=\{(\sigma_1, \dots, \sigma_p)\in \I_p : \prod_{j\in J_i}
\sigma_j =1,\ 1\leq i\leq r\} 
\ee
(a subgroup of $\I_p^+$).
\end{notation}

The groups $\I_\J^+$ generalize $\I_p^+$; we have $\I_\jtop^+=\I_p^+$.
Observe that an equivalent definition of  the group $\Gamma_\J^0$ defined in
Notation \ref{defkj} is $\Gamma_\J^0 = \{(\boldsig,\pi)\in\tsp^+ : \boldsig\in \I_\J^+,
\pi\in K_\J\}$. Thus, analogously to \eqref{ses2}, we have a short exact sequence
\bearray
\label{ses4}
1\to \I_\J^+ \to \Gamma_\J^0 \to  K_\J \to 1.
\eearray

Next, observe that the action \eqref{act1} of $\tsp^+$ on $M(p)$ induces, for each
$X\in\sympp$, an action of  $\tsp^+$ on $\Comp(\E_X)$, given by
\be\label{act3}
g\dotprod[(U,D)]:=[g\dotprod(U,D)].
\ee
\noi This leads us to:

\begin{prop}\label{prop:biject}
Let $X\in\sympp$. Then every $(U,D)\in\E_X$ determines a
  bijection between $\Comp(\E_X)$ and the set
  $\tsp^+/\Gamma^0_{\J_D}$.
\end{prop}

\pf Two elements $(U,D), (U',D')$ lie in the same component of $\E_X$
if and only if and only if $D'=D$ and $U'=UR$ for some $R\in G_D^0$.
Thus it is clear from \eqref{fibform4} that the action \eqref{act3} of
$\tsp^+$ on $\Comp(\E_X)$ is transitive.  Therefore for any $(U,D)\in
\E_X$, the map $\tsp^+\to \Comp(\E_X), g\mapsto g\dotprod[(U,D)]$,
induces a bijection $\tsp^+/\Stab([(U,D)])\to \Comp(\E_X)$, where
$\Stab([(U,D)])$ is the stabilizer of $[(U,D)]$ under the action
\eqref{act3}. But, as is easily checked,  $\Stab[(U,D)]$ is exactly the group
$\Gamma_{\J_D}^0$. \qedns

An important special case of Proposition \ref{prop:biject} is the case
in which
all eigenvalues of $X$ are distinct.  In this case, $\J_D=\jtop = \{
\{1\}, \{2\}, \dots, \{p\} \}$ and $\Gamma^0_{\J_D} =\{\id\}$. Thus the
action of $\tsp^+$ on $\Comp(\E_X)$ is free as well as transitive.
Furthermore $G_D^0=\{I\}$, so each connected component of $\E_X$ is
a single point; $\Comp(\E_X)=\E_X$. Thus $\E_X$ itself is an orbit of
$\tsp^+$, and any choice of $(U,D)\in \E_X$ yields a bijection
$\tsp^+\to \E_X$, $g\mapsto g\dotprod (U,D)$.

\begin{cor}\label{fibdescrip} Let $X\in \sympp$, $(U,D)\in \E_X$, and
let $k_1,\dots k_r$ be the parts of the partition $[\J_D]$ of $p$.
Then $\E_X$ is diffeomorphic to a disjoint union of
$2^{r-1}\frac{p!}{k_1! k_2! \dots k_r!}$ copies of $SO(k_1)\times
SO(k_2)\times \dots \times SO(k_r)$.
\end{cor}

\pf Let $\J=\J_D$. It is clear from \eqref{defcompud} that each
connected component of $\E_X$ is a submanifold of $M(p)$ diffeomorphic to
$G_D^0=G_\J^0$, which from \eqref{gj0iso} is isomorphic (hence
diffeomorphic) to $SO(k_1)\times SO(k_2)\times \dots \times
SO(k_r)$. From Proposition \ref{prop:biject}, the number of connected
components is $|\tsp^+/\Gamma^0_{\J_D}| = |\tsp^+|/
|\Gamma^0_{\J_D}|$.  As noted earlier, $|\tsp^+| =
2^{p-1}p!$, while from \eqref{ses4} we have $|\Gamma^0_{\J_D}| =
|\I_\J^+|\,|K_\J|$.  It is easily seen that $\I_\J^+$ is isomorphic to
$(\ztwo)^{p-r}$, and that $K_\J$ is isomorphic to $S_{k_1}\times
S_{k_2}\times \dots \times S_{k_r}$, and hence that $|K_\J|=k_1! k_2!
\dots k_r!$.  The result follows.  \qedns

\begin{remark}{\rm
An alternate, instructive route to Corollary
\ref{fibdescrip} is the following. (We merely sketch the ideas; the
reader may fill in the details.) For $\J\in\partpset$, define $ \Q_\J
= \{P_g R: g\in \tsp^+, R\in G_\J\}\subset \So(p)$.
Thus the set $\Q_\J$ is a
finite union of left-cosets of $G_\J$, each of which is diffeomorphic
to the compact submanifold $G_\J\subset \So(p)$.  If $X\in\sympp,$
$(U,D)\in \E_X$, and $\J=\J_D$, the map $\Q_\J\to M(p)$, $Q\mapsto
(UQ^{-1},QDQ^{-1})$, is an embedding with image $\E_X$. Hence $\E_X$
is a submanifold of $M(p)$ diffeomorphic to $\Q_\J$. But for any closed
subgroups $H_1, H_2$ of a compact Lie group $G$, the set
$H_1H_2:=\{h_1h_2: h_1\in H_1, h_2\in H_2\}\subset G$ is a submanifold
of $G$ and a principal $H_2$-bundle over $H_1/(H_1\intersect H_2)$,
with projection map given by $h_1h_2\mapsto h_1(H_1\intersect H_2)$.
Applying this to the case $H_1=\tsp^+, H_2=G_\J$, $G=\Sop$, we have
$H_1\intersect H_2=\Gamma_\J$, so $\Q_\J$ is a principal
$G_\J$-bundle over the {\em finite} set $\tsp^+/\Gamma_\J$. But the
natural map $\tsp^+/\Gamma_\J \to S_p/K_\J, \ g\Gamma_\J \mapsto
\proj_2(g)K_\J$ (where $\proj_2$ is as in \eqref{ses2}), is a bijection, so $\Q_\J$ may be viewed as a
principal $G_\J$-bundle over $S_p/K_\J$. The cardinality of this
base-space is $|S_p|/|K_\J|$, which is simply the multinomial
coefficient $\frac{p!}{k_1!  k_2! \dots k_r!}$ if $[\J]=(k_1,\dots,
k_r)\in \partp$. Thus $\E_X$ is diffeomorphic to $\frac{p!}{k_1!  k_2!
  \dots k_r!}$ copies of $G_\J$, and each copy of $G_\J$ is
diffeomorphic to $2^{r-1}$ copies of $SO(k_1)\times \dots \times
SO(k_r)$.  }
\end{remark}

\setcounter{equation}{0}
\setcounter{theorem}{0}
\renewcommand{\thesection}{Appendix B}

\section{Stratification of $\sympp$, $M(p)$, and related spaces}
\label{sect:strat}
\renewcommand{\thesection}{B}

We provide here a brief outline of the stratifications relevant to this paper.
For a more detailed discussion, see \cite[Section 2.7]{GJS2017eigen_EJS}.

 As noted in Section \ref{sec:notn_prelims}, $\Sop$ acts on $\sympp$ via
$(U,X)\mapsto UXU^T.$ As with any group-action, elements $X,Y\in\sympp$ are
said to have the same {\em orbit type} if their stabilizers are conjugate; in this
case the fibers $\E_X,\E_Y$ are diffeomorphic.  The orbit-type stratification of
any manifold under the action of a compact Lie group is
known to be a Whitney stratification (\cite[p. 21]{Gibsonetal1976}).

We use $\partpset$ to define stratifications of the spaces $\Dpp$ and $M(p)$,
and use $\partp$ to define stratifications of $\Dpp/S_p$ and $\sympp$. The
commutative diagram in Figure \ref{comdiag} indicates the relationships among
these spaces and label-sets. We define strata as the diagram suggests: for
$\J\in\partpset$ and $[\K]\in\partp$, (i) $\cald_\J := \lbl^{-1}(\J) \subset \Dpp$, (ii)
$\cald_{[\K]} := \Bar{\lbl}^{\,-1}([\K]) \subset \Dpp/S_p$, (iii) $\S_\J
:= \proj_2^{-1}(\cald_\J) = \Sop\times \cald_\J \subset M(p),$ and (iv)
$\S_{[\K]} :=\Bar{\proj_2}^{\,-1}(\cald_{[\K]})$.
The maps $\lbl, \Bar{\lbl}$ label elements of $\D^+(p),
\linebreak \D^+(p)/S_p$ by partitions of
the set $\{1,\dots,p\}$ and the integer $p$, respectively; $\proj_2: M(p)=
\So(p)\times\D^+(p)\to\D^+(p)$ is projection onto the second factor; and
$\Bar{\proj_2}$ is the map induced by $\proj_2$ on the indicated quotients.

For $X\in \S_{[\K]},$ we may call the partition $[\K]\in\partp$ the
{\em eigenvalue-multiplicity type} of $X$. The stratification  of $\sympp$ by
eigenvalue-multiplicity type is identical to the orbit-type stratification.

\begin{figure}[h]
\begin{diagram}
M(p)
& \rTo^{\ \sproj} &\D^+(p) & \rTo^{\slbl\ \ \ } &
\partpset\\
\dTo >F  &                &             \dTo >{\quo_1} &  &
\dTo >{\quo_2}\\
\sympp                  & \rTo^{\ \ \ \bsproj\ \ \ } & \D^+(p)/S_p &
\rTo^{\ \ \ \ \bslbl} & \partp\\
\end{diagram}
\caption{Commutative diagram defining the stratifications of $\sympp$ and related spaces.
}
\label{comdiag}
\end{figure}

In any stratified space, there is a natural partial ordering $\leq$ on the set of strata
${\cal T}_i$ defined by declaring ${\cal T}_1 \leq {\cal T}_2$ if ${\cal T}_1\subset
\Bar{{\cal T}_2}$.  Using this partial ordering of strata for the spaces in the left-hand
square in Figure \ref{comdiag}, it is easily checked that all the maps in Figure
\ref{comdiag} are either order-preserving themselves (in the case of $\quo_2$) or
induce order-preserving maps on the corresponding sets of strata (in the case of all
the other maps). In particular, each of the stratified spaces in the left-hand square in
Figure \ref{comdiag} has a top stratum and a bottom stratum.


\section*{References}

\begin{thebibliography}{10}
\expandafter\ifx\csname url\endcsname\relax
  \def\url#1{\texttt{#1}}\fi
\expandafter\ifx\csname urlprefix\endcsname\relax\def\urlprefix{URL }\fi
\expandafter\ifx\csname href\endcsname\relax
  \def\href#1#2{#2} \def\path#1{#1}\fi

\bibitem{Billera2001}
L.~J. Billera, S.~P. Holmes, K.~Vogtmann,
\href{http://dx.doi.org/10.1006/aama.2001.0759}
{Geometry of the space of phylogenetic trees}, {\em Adv. in Appl. Math.}
27~(4) (2001) 733--767.

\bibitem{CE1975}
J.~Cheeger, D.~G. Ebin, {\em Comparison Theorems in Riemannian Geometry}, North
  Holland/American Elsevier, Amsterdam, 1975.

\bibitem{Damon2013}
J.~Damon, J.~Marron, {Backwards principal
  component analysis and principal nested relations},
{\em J.   Math. Imaging and Vision} 50~(1) (2014), 107--114.

\bibitem{EAS1998}
A.~Edelman, T.~A. Arias, S.~T. Smith, The geometry of algorithms with
  orthogonality constraints, {\em SIAM J. Matrix Anal. Appl.} 20~(2) (1998) 303--353.

\bibitem{Gibsonetal1976}
C.~G. Gibson, K.~Wirthm{\"u}ller, A.~A. du~Plessis, E.~J.~N. Looijenga,
  {\em Topological Stability of Smooth Mappings}, Lecture Notes in Mathematics, Vol.
  552, Springer-Verlag, Berlin, 1976.

\bibitem{Golub1989}
G.~H. Golub, C.~F. {Van Loan}, {\em Matrix Computations}, 2nd edition, The Johns
  Hopkins University Press, 1989.

\bibitem{GJS2017eigen_EJS}
D.~Groisser, S.~Jung, A.~Schwartzman,  Geometric foundations for scaling-rotation
statistics on symmetric positive definite matrices: minimal smooth
scaling-rotation curves in low dimensions, {\em Electronic J. Stat.} {\bf 11~(1)},
1092--1159.

\bibitem{GJS2016metric}
D.~Groisser, S.~Jung, A.~Schwartzman, A scaling-rotation metric on the space of
  symmetric positive-definite matrices, in preparation.

\bibitem{Hotz2013}
T.~Hotz, S.~Huckemann, H.~Le, J.~S. Marron, J.~C. Mattingly, E.~Miller,
  J.~Nolen, M.~Owen, V.~Patrangenaru, S.~Skwerer,
  \href{http://dx.doi.org/10.1214/12-AAP899}{Sticky central limit theorems on
  open books}, {\em Ann. Appl. Prob.} 23~(6) (2013) 2238--2258.

\bibitem{JSG2015scarot}
S.~Jung, A.~Schwartzman, D.~Groisser, Scaling-rotation distance and
  interpolation of symmetric positive-definite matrices, {\em Siam J. Matrix Anal.
  Appl.}, 36~(3) (2015) 1180--1201.

\bibitem{Kendall1999}
D.~G. Kendall, D.~Barden, T.~K. Carne, H.~Le,
  \href{http://dx.doi.org/10.1002/9780470317006}{\em Shape and Shape Theory}, Wiley
  Series in Probability and Statistics, John Wiley \& Sons Ltd., Chichester,
  1999.

\bibitem{pahlings} H. Pahlings, Characterization of groups by their character tables,
{\em Comm. Alg.} 4~(2) (1976), 111--153.

\bibitem{samelson} H. Samelson, {\em Notes on Lie Algebras}, Van Nostrand Reinhold
Company, 1969.

\bibitem{Schwartzman2006}
A.~Schwartzman, Random ellipsoids and false discovery rates: statistics for
  diffusion tensor imaging data, Ph.D. thesis, Stanford University (2006).

\bibitem{Schwartzmanetal2008}
A.~Schwartzman, W.~F. Mascarenhas, J.~E. Taylor, Inference for eigenvalues and
  eigenvectors of {G}aussian symmetric matrices,
{\em Ann. Statist.}  36~(6) (2008) 2886--2919.

\bibitem{Wong1967}
Y.-C.~Wong, Differential geometry of {G}rassmann manifolds, {\em Proc. Nat. Acad.
  Sci. U.S.A.} 57 (1967) 589--594.

\end{thebibliography}
\end{document}